\documentclass[12pt]{article}

\makeatletter
\newfont{\footsc}{cmcsc10 at 8truept}
\newfont{\footbf}{cmbx10 at 8truept}
\newfont{\footrm}{cmr10 at 10truept}
\renewcommand{\ps@plain}{%
\renewcommand{\@oddfoot}{\footsc t.\ halverson and t.\ lewandowski, rsk insertion for set partitions, \today \hfil\footrm\thepage}}
  
\makeatother
\usepackage{amsmath,amsthm,amsfonts,amssymb}
\newtheorem{thm}{Theorem}[section]
\newtheorem{cor}[thm]{Corollary}

\newtheorem{prop}[thm]{Proposition}
\newtheorem{rem}[thm]{Remark}
\theoremstyle{definition}

\newtheorem{xmp}[thm]{Example}
\numberwithin{equation}{section} 
\newcommand{\ZZ}{\mathbb{Z}}
\newcommand{\CC}{\mathbb{C}}

\newcommand{\Res}{\mathrm{Res}}
\newcommand{\Ind}{\mathrm{Ind}}

\newcommand{\PT}{\mathcal{VT}}

\newcommand{\SYT}{\mathcal{SYT}}
\newcommand{\Card}{\mathrm{Card}}
\newcommand{\pn}{\mathrm{pn}}
\newcommand{\half}{\scriptstyle{\frac{1}{2}}}
\newcommand{\es}{{\scriptstyle{\emptyset}}}
\newcommand{\flip}{{\mathrm{flip}}}
\def\keyw#1{{\rm #1}}

\newcommand{\DI}{{\buildrel\rm{DI} \over \longrightarrow}}
\newcommand{\RS}{{\buildrel\rm{RSK} \over \longrightarrow}}
\newcommand{\jdt}{{\buildrel\rm{jdt} \over \longleftarrow}}

\input amssym.def
\input prepictex
\input pictex
\input postpictex

\begin{document}

\title {RSK Insertion for Set Partitions and \\ Diagram Algebras}
\author{Tom Halverson\thanks{Research supported in part by National Science Foundation Grant DMS0401098.}~ and  Tim Lewandowski${}^\ast$  \\ 
\small Department of Mathematics and Computer Science\\[-0.8ex]
\small Macalester College, Saint Paul, MN 55105 USA\\[-0.8ex]
%\small \texttt{halverson@macalester.edu} \qquad 
%\texttt{lewandowski@macalester.edu}
}

\maketitle

\centerline{\it In honor of Richard Stanley on his 60th birthday.}

\begin{abstract} We give combinatorial proofs of two identities from the representation
theory of the partition algebra $\CC A_k(n), n \ge 2k$.  The first is $n^k = \sum_\lambda
f^\lambda m_k^\lambda$, where the sum is over partitions $\lambda$ of $n$, $f^\lambda$
is the number of standard tableaux of shape $\lambda$, and $m_k^\lambda$ is the number of
``vacillating tableaux" of shape $\lambda$ and  length $2k$.  Our proof uses a combination of Robinson-Schensted-Knuth insertion and
jeu de taquin.  The second identity is $B(2k) = \sum_\lambda (m_k^\lambda)^2$, where $B(2k)$ is 
the number of set partitions of $\{1, \ldots, 2k\}$.   We show that this insertion restricts to work for the diagram algebras which appear as subalgebras of the partition algebra: the Brauer, Temperley-Lieb, planar partition, rook monoid, planar rook monoid, and symmetric group algebras.
\end{abstract}

\section{Introduction}
Two fundamental identities in the representation theory of
the symmetric group $S_k$ are
\begin{equation}
{\rm (a)} \qquad n^k = \sum_{\lambda \vdash k \atop \ell(\lambda) \le n}  f^\lambda d_\lambda,
\qquad\hbox{ and }\qquad
{\rm (b)} \qquad k! =\sum_{\lambda \vdash k}  (f^\lambda)^2, \qquad
\label{eq:Skids} 
\end{equation}
where $\lambda$ varies over partitions of the integer $k$ of length $\ell(\lambda) \le n$,
$f^\lambda$ is the number of standard Young tableaux of shape $\lambda$, 
and $d_\lambda$ is the number of column strict tableaux of shape $\lambda$
with entries from $\{1, \ldots, n\}$.  The Robinson-Schensted-Knuth (RSK) insertion algorithm provides a  a bijection between sequences $(i_1, \ldots, i_k), 1 \le i_j \le n,$ and pairs $(P_\lambda,Q_\lambda)$ consisting of a standard Young tableau $P_\lambda$ of shape $\lambda$ and a column strict tableau $Q_\lambda$ of shape $\lambda$, thus providing a combinatorial proof of  (\ref{eq:Skids}.a).  If we restrict  $i_1, \ldots, i_k$ to be a permutation of $1, \ldots, k$, then $Q_\lambda$ is  a standard tableau and we have a proof of (\ref{eq:Skids}.b).

Identity (\ref{eq:Skids}.b) comes from the decomposition of the
group algebra $\CC[S_k]$ into irreducible $S_k$-modules $V^\lambda, \lambda\vdash k$,
where $\dim(V^\lambda) = f^\lambda$ and the multiplicity of $V^\lambda$ in $\CC[S_k]$
is also $f^\lambda$. Identity (\ref{eq:Skids}.a)
comes from the Schur-Weyl duality between $S_k$ and the general linear
group $GL_n(\CC)$ on the $k$-fold tensor product $V^{\otimes k}$ of the
fundamental representation $V$ of $GL_n(\CC)$.  There is an action of $S_k$
on $V^{\otimes k}$ by tensor place permutations, and via this action, $\CC[S_k]$ is isomorphic
to the centralizer algebra $End_{GL_n(\CC)}(V^{\otimes k})$. As a bimodule for 
$S_k \times GL_n(\CC)$,
\begin{equation}
V^{\otimes k} \cong  \bigoplus_{\lambda \vdash k}\  S^\lambda \otimes V^\lambda,
\label{eq:Skdecomp}
\end{equation}
where $S^\lambda$ is an irreducible $S_k$-module of dimension $f^\lambda$ and
$V^\lambda$ is an irreducible $GL_r(\CC)$-module of dimension $d_\lambda$. 
We get  (\ref{eq:Skids}.a) by computing dimensions on each side of (\ref{eq:Skdecomp}).

R.\ Brauer \cite{Br} defined an algebra $\CC B_k(n)$, which is isomorphic to the centralizer algebra of the orthogonal group $O_n(\CC) \subseteq GL_n(\CC)$, when $n \ge 2k$, i.e., 
$
\CC B_k(n) \cong  End_{O_n(\CC)}(V^{\otimes k}).
$
The dimension of the Brauer algebra is $(2k-1)!! = (2k-1)(2k-3) \cdots 3 \cdot 1$.
Since $O_n(\CC) \subseteq GL_n(\CC)$, their centralizers satisfy $\CC B_k(n) \supseteq \CC[S_k]$.
Berele \cite{Be} generalized the RSK correspondence to give combinatorial
proof of the $\CC B_k(n)$-analog of (\ref{eq:Skids}.a).  Sundaram \cite{Sun} (see also \cite{Ter}) gave
a combinatorial proof of the $\CC B_k(n)$-analog  of (\ref{eq:Skids}.a).

We now take this restriction  further to $S_{n-1} \subseteq S_n \subseteq O_n(\CC) \subseteq GL_n(\CC)$, where $S_n$ is viewed as the subgroup of permutation matrices in $GL_n(\CC)$ and $S_{n-1} \subseteq S_n$ corresponds to the permutations that fix $n$. Under this restriction, $V$ is the permutation representation of $S_n$, and when $n \ge 2k$, the centralizer algebras are the partition algebras,
$$
\CC A_k(n) \cong End_{S_n}(V^{\otimes k}) \quad\hbox{and}\quad
\CC A_{k+{1\over2}}(n) \cong End_{S_{n-1}}(V^{\otimes k}).
$$
The partition algebra $\CC A_k(n)$ first appeared independently 
in the work of  Martin \cite{Mar1,Mar2,Mar3}  and Jones \cite{Jo} arising from applications in statistical
mechanics. See \cite{HR2} for a survey paper on partition algebras.

For $k \in {1 \over 2} \ZZ_{>0}$ and $n \ge 2 k$,  $\CC A_k(n)$ is  a semisimple algebra over $\CC$ with a basis  indexed by the set partitions of $\{1, \ldots, 2k\}$.  Thus, the dimension of $\CC A_k(n)$  is the $2k$th Bell number $B(2k)$.  Define
 \begin{equation}
\Lambda_n^k = \left\{ \lambda \vdash n\, \Big| \, |\lambda| - \lambda_1 \le k\,\right\}
\label{eq:irreducibles}
\end{equation}
(these are partitions of $n$ with at most $k$ boxes below the first row of their Young diagram).
Then for $k \in \ZZ_{>0}$, the irreducible representations of $\CC A_k(n)$ are
indexed by partitions in the set $\Lambda_n^k$, and the irreducible  representations of $\CC A_{k+{1 \over 2}}(n)$ are indexed by partitions in the set $\Lambda_{n-1}^k$.

Using the Schur-Weyl duality between $S_n$ and $\CC A_k(n)$  we get the identity
\begin{equation}
n^k = \sum_{\lambda\in \Lambda_n^k} f^\lambda m_k^\lambda, 
\label{eq:Pkid1}
\end{equation}
where $m_k^\lambda$ is the number of vacillating tableaux of shape $\lambda$ and  length $2k$ (defined in
Section 2.2), which are sequences integer partitions in the Bratteli diagram of $\CC A_k(n)$. Identity
(\ref{eq:Pkid1})
is the partition algebra analog of (\ref{eq:Skids}.a). In Section 3,  we prove
(\ref{eq:Pkid1}) using RSK column insertion and jeu de taquin.
Decomposing $\CC A_k(n)$ as a bimodule for $\CC A_k(n) \otimes \CC A_k(n)$ gives
\begin{equation}
B(2k) = \sum_\lambda (m_k^\lambda)^2,
\label{eq:Pn2}
\end{equation}
which is the $\CC A_k(n)$-analog of (\ref{eq:Skids}.b). In Section 4, we give a bijective proof of 
(\ref{eq:Pn2}) that contains as a special cases the RSK algorithms for $\CC S_k$ and $\CC B_k(n)$. 

Martin
and Rollet \cite{MR} have given a {\it different\/} combinatorial proof of the second identity
(\ref{eq:Pn2}). Their bijection has the elegant property that pairs of  paths in the
difference between the Bratteli diagrams of $\CC A_k(\ell)$ and $\CC A_k(\ell -1)$ 
are in exact correspondence with the set partitions of $\{1, \ldots, 2k\}$
into $\ell$ parts.  The advantages of the correspondence in this paper are:
\begin{enumerate}\addtolength{\itemsep}{-0.3\baselineskip}
\item Our algorithm for (\ref{eq:Pn2}) contains, as special cases, the known RSK  correspondences for a number of  diagram algebras which appear as subalgebras of $\CC A_k(n)$:
\begin{enumerate}\addtolength{\itemsep}{-0.3\baselineskip}
\item The group algebra of the symmetric group $\CC S_k$,
\item The Brauer algebra $\CC B_k(n)$,
\item The Temperley-Lieb algebra $\CC T_k(n)$,
\item The planar partition algebra $\CC P_k(n)$.
\item The rook monoid algebra $\CC R_n$ and the planar rook monoid algebra $\CC P\!R_n$. 
\end{enumerate}
Thus we obtain combinatorial proofs of the analog of (\ref{eq:Pn2}) for each of these algebras (see equations (5.1), (5.2), (5.3), and (5.4)).

\item  Using Fomin growth diagrams we show  that our algorithm for (\ref{eq:Pn2}) is symmetric in the sense that if $d \rightarrow (P,Q)$ then $\flip(d) \mapsto (Q,P)$ where $\flip(d)$ is the diagram $d$ flipped over its horizontal axis.  This is the generalization of the property for the symmetric group that if $\pi \rightarrow (P,Q)$ then $\pi^{-1} \rightarrow (Q,P)$.  As a consequence, we show that the number of symmetric diagrams equals the sum of the dimensions of the irreducible representations (for each of the diagram algebras mentioned in item 1 above).
\item  Our algorithms for the bijections in (\ref{eq:Pkid1}) and (\ref{eq:Pn2}) each use iterations of RSK insertion and jeu de taquin (see equations (3.4) and (4.4)).

\end{enumerate}

The main idea for the algorithm in this paper came from a bijection of R.\ Stanley between fixed point free involutions in the symmetric group $S_{2k}$ and Brauer diagrams. This led to the insertion scheme of Sundaram in \cite{Sun} for the Brauer algebra. After we distributed a preliminary version  of this paper, R. Stanley and colleagues independently came out with the  paper \cite{CDDSY}, which studies the  crossing and nesting properties of an extended version of the insertion used in this paper for the bijection in (1.5).  We have adopted the term ``vacillating tableaux" from \cite{CDDSY}, and we use the crossing property of  \cite{CDDSY} to show that our algorithm restricts appropriately to the planar partition algebra.

%%%%%%%%%%%%%%%%%%%%%%%%%%%%%%%%%%%%%%%%%%%%%%%

\section{The Partition Algebra and Vacillating Tableaux}

\bigskip
For $k \in \ZZ_{> 0}$, let
\begin{eqnarray} 
A_k & =& \big\{ \hbox{set partitions of
$\{1,2,\ldots,k,1',2',\ldots, k'\}$}\big\},
\qquad\hbox{and}\\
A_{k+{1 \over 2}} & =&
\left\{ d \in \Pi_{k+1} ~ \big| ~
\hbox{$(k+1)$ and $(k+1)'$ are in the same block}
\right\}.
\end{eqnarray}
The propagating number of $d \in  A_k$ is
\begin{equation}\label{def:propagating}
\pn(d) = 
\left(
\begin{array}{l}
\hbox{the number of blocks in $d$ that contain both an element}\\
\hbox{of $\{1, 2, \ldots, k\}$ and an element of $\{1',2',\ldots,k'\}$}
\end{array}\right).
\end{equation}
For convenience, represent a set partition $d\in A_k$
by a graph with $k$ vertices in the top row,
labeled $1, \ldots,k$,  and $k$ vertices
in the bottom row, labeled $1', \ldots, k'$,
with vertex $i$ and vertex $j$ connected by a path if $i$
and $j$ are in the same block of the set partition $d$.
For example,
$$
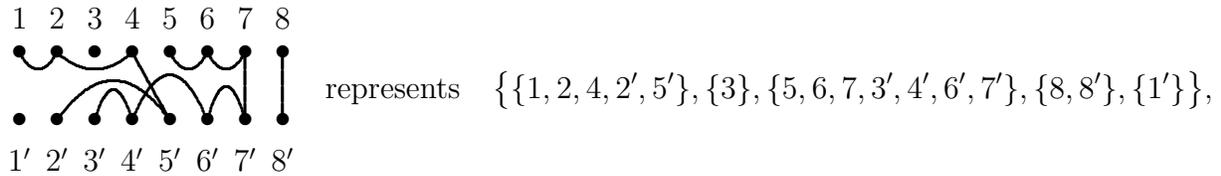
{\beginpicture
\setcoordinatesystem units <0.5cm,0.3cm> % sets scale
\setplotarea x from 0 to 7, y from 0 to 3    % sets plot size up
\linethickness=0.5pt                        % sets line thickness
\put{1} at 0 3.5
\put{2} at 1 3.5
\put{3} at 2 3.5
\put{4} at 3 3.5
\put{5} at 4 3.5
\put{6} at 5 3.5
\put{7} at 6 3.5
\put{8} at 7 3.5
\put{$1'$} at  0 -2.75
\put{$2'$} at 1 -2.75
\put{$3'$} at 2 -2.75
\put{$4'$} at 3 -2.75
\put{$5'$} at 4 -2.75
\put{$6'$} at 5 -2.75
\put{$7'$} at 6 -2.75
\put{$8'$} at 7 -2.75
\put{$\bullet$} at 0 -1 \put{$\bullet$} at 0 2
\put{$\bullet$} at 1 -1 \put{$\bullet$} at 1 2
\put{$\bullet$} at 2 -1 \put{$\bullet$} at 2 2
\put{$\bullet$} at 3 -1 \put{$\bullet$} at 3 2
\put{$\bullet$} at 4 -1 \put{$\bullet$} at 4 2
\put{$\bullet$} at 5 -1 \put{$\bullet$} at 5 2
\put{$\bullet$} at 6 -1 \put{$\bullet$} at 6 2
\put{$\bullet$} at 7 -1 \put{$\bullet$} at 7 2
%\plot 1 -1 1 2  /
\plot 3 2 4 -1 /
\plot 6 -1 6 2  /
\plot 7 -1 7 2  /
\setquadratic
\plot 0  2  .5 1.25 1 2 /
\plot 4  2 4.5 1.25 5 2 /
\plot 5  2 5.5 1.25 6 2 /
%\plot 0  2  .5 .25 1 -1 /
\plot 1 -1 2.5 .75 4 -1 /
\plot 1  2  2  1.25 3 2 /
\plot 2 -1 2.5 0.35 3 -1 /
\plot 3 -1 4 1 5 -1 /
\plot 5 -1 5.5 .5 6 -1 /
\endpicture}
\quad\hbox{represents}\quad
\big\{ \{ 1, 2, 4, 2',5'\},\{3\},
\{5, 6, 7, 3', 4', 6', 7'\},\{8, 8'\}, \{1'\}\big\},
$$
and has propagating number 3.
The graph representing $d$ is not unique.

Define the composition $d_1\circ d_2$ of partition diagrams
$d_1,d_2\in A_k$ to be the set partition $d_1 \circ d_2\in A_k$
obtained by placing $d_1$ above $d_2$, identifying the bottom
dots of $d_1$ with the top dots of $d_2$, and removing any connected 
components that live entirely in the middle row.  For example, 
$$
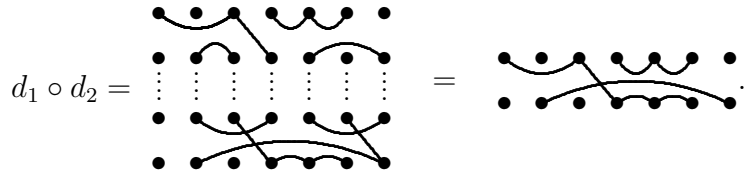
{\beginpicture
\setcoordinatesystem units <0.5cm,0.2cm> % sets scale
\setplotarea x from 0 to 7, y from 0 to 3    % sets plot size up
\linethickness=0.5pt
\put{$d_1 \circ d_2 = $\qquad} at -.5 0
\put{$\bullet$} at 1 2 \put{$\bullet$} at 1 5
\put{$\bullet$} at 2 2 \put{$\bullet$} at 2 5
\put{$\bullet$} at 3 2 \put{$\bullet$} at 3 5
\put{$\bullet$} at 4 2 \put{$\bullet$} at 4 5
\put{$\bullet$} at 5 2 \put{$\bullet$} at 5 5
\put{$\bullet$} at 6 2 \put{$\bullet$} at 6 5
\put{$\bullet$} at 7 2 \put{$\bullet$} at 7 5

\put{$\bullet$} at 1 -2 \put{$\bullet$} at 1 -5
\put{$\bullet$} at 2 -2 \put{$\bullet$} at 2 -5
\put{$\bullet$} at 3 -2 \put{$\bullet$} at 3 -5
\put{$\bullet$} at 4 -2 \put{$\bullet$} at 4 -5
\put{$\bullet$} at 5 -2 \put{$\bullet$} at 5 -5
\put{$\bullet$} at 6 -2 \put{$\bullet$} at 6 -5
\put{$\bullet$} at 7 -2 \put{$\bullet$} at 7 -5

\setquadratic
\plot 1 5 2 4.0 3 5 /
\plot 1 5 2 4.0 3 5 /
\plot 4 5 4.5 4.0 5 5 /
\plot 5 5 5.5 4.0 6 5 /
\plot 2 2 2.5 3 3 2 /
\plot 5 2 6 3 7 2 /
\plot 2 -2 3 -3 4 -2 /
\plot 5 -2 6 -3 7 -2 /
\plot 2 -5 4.5 -3.5 7 -5 /
\plot 4 -5 4.5 -4.5 5 -5 /
\plot 5 -5 5.5 -4.5 6 -5 /
\setlinear
\plot 3 5 4 2 /
\plot 3 -2 4 -5 /
\plot 6 -2 7 -5 /
\setdashpattern <.01cm,.1cm>
\plot 1 1 1 -1 /
\plot 2 1 2 -1 /
\plot 3 1 3 -1 /
\plot 4 1 4 -1 /
\plot 5 1 5 -1 /
\plot 6 1 6 -1 /
\plot 7 1 7 -1 /
\endpicture}
\quad = \quad {\beginpicture
\setcoordinatesystem units <0.5cm,0.2cm> % sets scale
\setplotarea x from 1 to 7, y from 0 to 3    % sets plot size up
\linethickness=0.5pt                        % sets line thickness
\put{$\bullet$} at 1 -1 \put{$\bullet$} at 1 2
\put{$\bullet$} at 2 -1 \put{$\bullet$} at 2 2
\put{$\bullet$} at 3 -1 \put{$\bullet$} at 3 2
\put{$\bullet$} at 4 -1 \put{$\bullet$} at 4 2
\put{$\bullet$} at 5 -1 \put{$\bullet$} at 5 2
\put{$\bullet$} at 6 -1 \put{$\bullet$} at 6 2
\put{$\bullet$} at 7 -1 \put{$\bullet$} at 7 2
\plot 3 2 4 -1 /
\setquadratic
\plot 1 2 2 1.0 3 2 /
\plot 4 2 4.5 1.0 5 2 /
\plot 5 2 5.5 1.0 6 2 /
\plot 2 -1 4.5 .5 7 -1 /
\plot 4 -1 4.5 -.5 5 -1 /
\plot 5 -1 5.5 -.5 6 -1 /
\endpicture}.
$$
Diagram multiplication makes  $A_k$ into an
associative monoid with identity,
$1 {\beginpicture
\setcoordinatesystem units <0.35cm,0.12cm> % sets scale
\setplotarea x from 1 to 4, y from -1 to 2   % sets plot size up
\linethickness=0.5pt                        % sets line thickness
\put{$\bullet$} at 1 -1 \put{$\bullet$} at 1 2 \plot 1 -1 1 2 /
\put{$\bullet$} at 2 -1 \put{$\bullet$} at 2 2 \plot 2 -1 2 2 /
\put{$\cdots$} at 3 .5
\put{$\bullet$} at 4 -1  \put{$\bullet$} at 4 2
\plot 4 -1 4 2 /
\endpicture},
$
and the propagating number satisfies $\pn(d_1 \circ d_2) \le \min(\pn(d_1),\pn(d_2)). $

For $k \in {1 \over 2} \ZZ_{> 0}$ and $n \in \CC$, the partition
algebra $\CC A_k(n)= \CC\hbox{span-}\{ d\in A_k\}$ is an associative algebra over 
$\CC$ with basis $A_k$.  Multiplication in  $\CC A_k(n)$ is defined by 
$$
d_1d_2 = n^\ell (d_1 \circ d_2),
$$
where  $\ell$ is the number of blocks removed from the the middle row when constructing
the composition $d_1 \circ d_2$.  In the example above $d_1 d_2 = n^2 d_1 \circ d_2$.

 For each $k \in {1 \over 2} \ZZ_{>0}$,
the following are submonoids of the partition monoid $A_k$:
\begin{eqnarray*}
S_k & = &  \{ d \in A_k \ | \ \pn(d) = k \},  \qquad\qquad
I_t   =   \{ d \in A_k \ | \ \pn(d) \le  t \},  \quad 0 < t \le k, \\
B_k & =&  \{ d\in A_k \ |\ \hbox{all blocks of $d$ have size 2}\}, \\
R_k &=& \left\{d\in A_k \ \bigg|\ 
\begin{array}{l}
\hbox{all blocks of $d$ have at most one vertex in $\{1, \ldots k\}$} \\
\hbox{and at most one vertex in $\{1', \ldots k'\}$} \\
\end{array}
\right\}.
\end{eqnarray*}
Here, $S_k$ is the symmetric group, $B_k$ is the Brauer monoid, and $R_k$ is the rook monoid.
A set partition is {\it planar} [Jo] if it can be represented
as a graph without edge crossings inside of the rectangle formed
by its vertices. The following are planar submonoids,
$$
P_k  =  \{ d \in A_k \ | \ \hbox{$d$ is planar} \}, \quad
T_k =  B_k \cap P_k, \quad
P\!R_k = R_k \cap P_k , \quad 1 = S_k \cap P_k.
$$
Here, $T_k$ is the Temperley-Lieb monoid.
Examples of diagrams in the various submonoids are:
$$
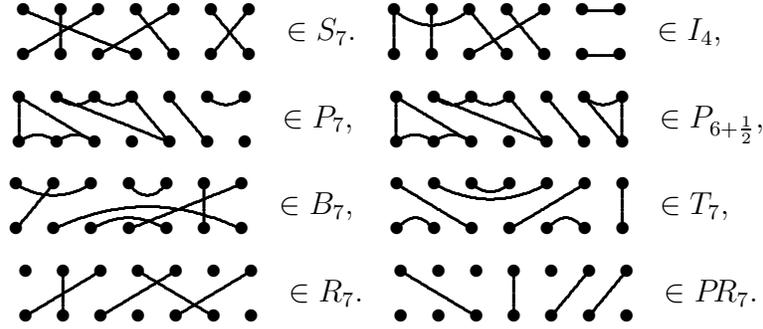
{\beginpicture
\setcoordinatesystem units <0.5cm,0.2cm> % sets scale
\setplotarea x from 0 to 7, y from 0 to 3    % sets plot size up
\linethickness=0.5pt                        % sets line thickness
\put{$\bullet$} at 1 -1 \put{$\bullet$} at 1 2
\put{$\bullet$} at 2 -1 \put{$\bullet$} at 2 2
\put{$\bullet$} at 3 -1 \put{$\bullet$} at 3 2
\put{$\bullet$} at 4 -1 \put{$\bullet$} at 4 2
\put{$\bullet$} at 5 -1 \put{$\bullet$} at 5 2
\put{$\bullet$} at 6 -1 \put{$\bullet$} at 6 2
\put{$\bullet$} at 7 -1 \put{$\bullet$} at 7 2
\plot 1 2 4 -1 /
\plot 2 2 2 -1 /
\plot 3 2 1 -1 /
\plot 4 2 5 -1 /
\plot 5 2 3 -1 /
\plot 6 2 7 -1 /
\plot 7 2 6 -1 /
\endpicture}
\ \ \in  S_7.
{\beginpicture
\setcoordinatesystem units <0.5cm,0.2cm> % sets scale
\setplotarea x from 0 to 7, y from 0 to 3    % sets plot size up
\linethickness=0.5pt                        % sets line thickness
\put{$\bullet$} at 1 -1 \put{$\bullet$} at 1 2
\put{$\bullet$} at 2 -1 \put{$\bullet$} at 2 2
\put{$\bullet$} at 3 -1 \put{$\bullet$} at 3 2
\put{$\bullet$} at 4 -1 \put{$\bullet$} at 4 2
\put{$\bullet$} at 5 -1 \put{$\bullet$} at 5 2
\put{$\bullet$} at 6 -1 \put{$\bullet$} at 6 2
\put{$\bullet$} at 7 -1 \put{$\bullet$} at 7 2
\plot 3 2 4 -1 /
\plot 2 2 2 -1 /
\plot 1 2 1 -1 /
\plot 4 2 5 -1 /
\plot 5 2 3 -1 /
\plot 6 -1 7 -1 /
\plot 7 2 6 2 /
\setquadratic
\plot 1 2 2 1 3 2 /
\endpicture}
\ \ \in  I_4,  \hskip.2truein
$$
$$
%%%%%%%
{\beginpicture
\setcoordinatesystem units <0.5cm,0.2cm> % sets scale
\setplotarea x from 0 to 7, y from 0 to 3    % sets plot size up
\linethickness=0.5pt                        % sets line thickness
\put{$\bullet$} at 1 -1 \put{$\bullet$} at 1 2
\put{$\bullet$} at 2 -1 \put{$\bullet$} at 2 2
\put{$\bullet$} at 3 -1 \put{$\bullet$} at 3 2
\put{$\bullet$} at 4 -1 \put{$\bullet$} at 4 2
\put{$\bullet$} at 5 -1 \put{$\bullet$} at 5 2
\put{$\bullet$} at 6 -1 \put{$\bullet$} at 6 2
\put{$\bullet$} at 7 -1 \put{$\bullet$} at 7 2
\plot 1 2 3 -1 /
\plot 5 2 6 -1 /
\plot 2 2 5 -1 /
\plot 4 2 5 -1 /
\plot 1 2 1 -1 /
\setquadratic
\plot 2 2 2.5 1.5 3 2 /
\plot 3 2 3.5 1.5 4 2 /
\plot 6 2 6.5 1.5 7 2 /
\plot 1 -1 1.5 -.5 2 -1 /
\plot 2 -1 2.5 -.5 3 -1 /
\endpicture}
\ \ \in P_{7}, 
%%%%%%%
{\beginpicture
\setcoordinatesystem units <0.5cm,0.2cm> % sets scale
\setplotarea x from 0 to 7, y from 0 to 3    % sets plot size up
\linethickness=0.5pt                        % sets line thickness
\put{$\bullet$} at 1 -1 \put{$\bullet$} at 1 2
\put{$\bullet$} at 2 -1 \put{$\bullet$} at 2 2
\put{$\bullet$} at 3 -1 \put{$\bullet$} at 3 2
\put{$\bullet$} at 4 -1 \put{$\bullet$} at 4 2
\put{$\bullet$} at 5 -1 \put{$\bullet$} at 5 2
\put{$\bullet$} at 6 -1 \put{$\bullet$} at 6 2
\put{$\bullet$} at 7 -1 \put{$\bullet$} at 7 2
\plot 1 2 3 -1 /
\plot 7 2 7 -1 /
\plot 6 2 7 -1 /
\plot 5 2 6 -1 /
\plot 2 2 5 -1 /
\plot 4 2 5 -1 /
\plot 1 2 1 -1 /
\setquadratic
\plot 2 2 2.5 1.5 3 2 /
\plot 3 2 3.5 1.5 4 2 /
\plot 6 2 6.5 1.5 7 2 /
\plot 1 -1 1.5 -.5 2 -1 /
\plot 2 -1 2.5 -.5 3 -1 /
\endpicture}
\ \ \in P_{6+{1 \over 2}}, 
$$
$$
%%%%%%%
{\beginpicture
\setcoordinatesystem units <0.5cm,0.2cm> % sets scale
\setplotarea x from 0 to 7, y from 0 to 3    % sets plot size up
\linethickness=0.5pt                        % sets line thickness
\put{$\bullet$} at 1 -1 \put{$\bullet$} at 1 2
\put{$\bullet$} at 2 -1 \put{$\bullet$} at 2 2
\put{$\bullet$} at 3 -1 \put{$\bullet$} at 3 2
\put{$\bullet$} at 4 -1 \put{$\bullet$} at 4 2
\put{$\bullet$} at 5 -1 \put{$\bullet$} at 5 2
\put{$\bullet$} at 6 -1 \put{$\bullet$} at 6 2
\put{$\bullet$} at 7 -1 \put{$\bullet$} at 7 2
\plot 7 2 4 -1 /
\plot 6 2 6 -1 /
\plot 2 2 1 -1 /
\setquadratic
\plot 1 2 2 1.25 3 2 /
\plot 4 2 4.5 1.25 5 2 /
\plot 2 -1 4.5 0.5 7 -1 /
\plot 3 -1 4 -.25 5 -1 /
\endpicture}
\ \ \in  B_7, 
%%%%%%%
{\beginpicture
\setcoordinatesystem units <0.5cm,0.2cm> % sets scale
\setplotarea x from 0 to 7, y from 0 to 3    % sets plot size up
\linethickness=0.5pt                        % sets line thickness
\put{$\bullet$} at 1 -1 \put{$\bullet$} at 1 2
\put{$\bullet$} at 2 -1 \put{$\bullet$} at 2 2
\put{$\bullet$} at 3 -1 \put{$\bullet$} at 3 2
\put{$\bullet$} at 4 -1 \put{$\bullet$} at 4 2
\put{$\bullet$} at 5 -1 \put{$\bullet$} at 5 2
\put{$\bullet$} at 6 -1 \put{$\bullet$} at 6 2
\put{$\bullet$} at 7 -1 \put{$\bullet$} at 7 2
\plot 1 2 3 -1 /
\plot 6 2 4 -1 /
\plot 7 2 7 -1 /
\setquadratic
\plot 3 2 3.5 1.5 4 2 /
\plot 2 2 3.5 1.0 5 2 /
\plot 1 -1 1.5 -.25 2 -1 /
\plot 5 -1 5.5 -.25 6 -1 /
\endpicture}
\ \ \in T_7, \hskip.2truein
$$
$$
%%%%%%%
{\beginpicture
\setcoordinatesystem units <0.5cm,0.2cm> % sets scale
\setplotarea x from 0 to 7, y from 0 to 3    % sets plot size up
\linethickness=0.5pt                        % sets line thickness
\put{$\bullet$} at 1 -1 \put{$\bullet$} at 1 2
\put{$\bullet$} at 2 -1 \put{$\bullet$} at 2 2
\put{$\bullet$} at 3 -1 \put{$\bullet$} at 3 2
\put{$\bullet$} at 4 -1 \put{$\bullet$} at 4 2
\put{$\bullet$} at 5 -1 \put{$\bullet$} at 5 2
\put{$\bullet$} at 6 -1 \put{$\bullet$} at 6 2
\put{$\bullet$} at 7 -1 \put{$\bullet$} at 7 2
\plot 2 2 2 -1 /
\plot 3 2 1 -1 /
\plot 4 2 6 -1 /
\plot 5 2 3 -1 /
\plot 7 2 5 -1 /
\endpicture}
\ \ \in  R_7.
%%%%%%%
{\beginpicture
\setcoordinatesystem units <0.5cm,0.2cm> % sets scale
\setplotarea x from 0 to 7, y from 0 to 3    % sets plot size up
\linethickness=0.5pt                        % sets line thickness
\put{$\bullet$} at 1 -1 \put{$\bullet$} at 1 2
\put{$\bullet$} at 2 -1 \put{$\bullet$} at 2 2
\put{$\bullet$} at 3 -1 \put{$\bullet$} at 3 2
\put{$\bullet$} at 4 -1 \put{$\bullet$} at 4 2
\put{$\bullet$} at 5 -1 \put{$\bullet$} at 5 2
\put{$\bullet$} at 6 -1 \put{$\bullet$} at 6 2
\put{$\bullet$} at 7 -1 \put{$\bullet$} at 7 2
\plot 1 2 3 -1 /
\plot 4 2 4 -1 /
\plot 6 2 5 -1 /
\plot 7 2 6 -1 /
\endpicture}
\ \ \in  P\!R_7.
$$

For each  monoid, we make an associative algebra in the same way that we construct the partition algebra $\CC A_k(n)$ from the partition monoid $A_k(n)$.  For example, we obtain the
the Brauer algebra $\CC B_k(n)$, the Temperley-Lieb algebra $\CC T_k(n)$, the group algebra of the symmetric group $\CC S_k$ in this way. Multiplication in the rook monoid algebra $\CC R_k$ is done without the coefficient $n^\ell$ (see \cite{Ha}).

For $\ell\in \ZZ_{>0}$, the Bell number  $B(\ell)$ is the number of set partitions of $\{1, 2, \ldots, \ell\}$,
the  Catalan number is
$
C(\ell) = {1 \over \ell + 1} {2 \ell \choose \ell}
= {2 \ell \choose \ell} - {2 \ell \choose \ell+ 1},
$
and 
$
(2 \ell)!! = (2\ell-1)\cdot (2 \ell -3) \cdots 5 \cdot 3 \cdot 1.
$
These have  generating functions (see \cite[1.24f, and 6.2]{Sta}),
$$
\sum_{\ell \ge 0} B(\ell) {z^\ell \over \ell !}  = \exp(e^z - 1),\
\qquad\qquad
\sum_{\ell \ge 0} C(\ell-1) z^\ell = {1 - \sqrt{1 - 4z}\over 2z},
$$
and
$$
\sum_{\ell \ge 0} (2 (\ell-1))!! {z^\ell \over \ell !} = {1 - \sqrt{1 - 2 z}\over z}.
$$
For $k \in {1 \over 2} \ZZ_{>0}$, the monoids have cardinality
$$
\Card(A_k) = B(2 k), \qquad
\Card(P_k) = \Card(T_{2k}) = C(2k), \qquad\hbox{for $k \in {1 \over 2} \ZZ_{>0}$}
$$
and 
$$
\begin{array}{lll}
\Card(S_k) =  k!,  & \qquad &  \Card(B_k) =  (2k)!!, \\
\Card(R_k) = \displaystyle{\sum_{\ell = 0}^k {k \choose \ell}^2 \ell!}
& \qquad & \Card(P\!R_k) = \displaystyle{{2 k \choose k}}, \\
\end{array}  \qquad\hbox{for $k \in \ZZ_{>0}$.}
$$

\subsection{Schur-Weyl Duality Between $S_n$ and $\CC A_k(n)$}

The irreducible representations of $S_n$ are indexed by integer partitions of $n$. If $\lambda (\lambda_1, \ldots,
\lambda_\ell) \in (\ZZ^{\ge 0})^\ell$ with $\lambda_1 \ge \cdots \ge \lambda_\ell$ and $\lambda_1 + \cdots +
\lambda_\ell = n$, then $\lambda$ is a partition of $n$, denoted $\lambda\vdash n$.  If $\lambda \vdash n$,
then we write $|\lambda| = n$. If $\lambda = (\lambda_1, \cdots, \lambda_\ell)$ and $\mu = (\mu_1, \ldots, \mu_\ell)$
are partitions such that $\mu_i \le \lambda_i$ for each $i$, then we say that $\mu \subseteq \lambda$,
and $\lambda/\mu$ is the skew shape given by deleting the boxes of $\mu$ from the Young diagram of $\lambda$.

Let $V$ be the $n$-dimensional permutation representation of the symmetric
group $S_n$. If we view $S_{n-1} \subseteq S_n$ as the subgroup of permutations that fix $n$, then $V$ is isomorphic to the left coset representation $\CC[S_n/S_{n-1}]$.  Let $V^{\otimes k}$ be the $k$-fold tensor product representation of $V$, and let $V^{\otimes 0} = \CC$.
From the ``tensor identity"  (see for example \cite{HR2}), we have the following restriction-induction rule
rule for $i \ge 0$,
\[
V^{\otimes (i+1)} \cong  V^{\otimes i} \otimes \CC[S_n/S_{n-1}] 
\cong
\Ind_{S_{n-1}}^{S_n}(\Res^{S_n}_{S_{n-1}}(V^{\otimes i})).
\]
Thus, $V^{\otimes k}$ is obtained from $k$ iterations of  restricting to $S_{n-1}$ and inducing back to $S_n$.

Let $V^\lambda$ denote the irreducible representations of $S_n$ indexed by $\lambda \vdash n$. 
The restriction and induction rules for $S_{n-1} \subseteq S_n$ are given by
\begin{eqnarray}
\Res^{S_n}_{S_{n-1}}(V^\lambda) &\cong& \bigoplus_{\mu \vdash (n-1), \, \mu \subseteq \lambda} V^\mu,
\qquad\hbox{for $\lambda \vdash n$}
 \label{eq:res}\\
\Ind^{S_n}_{S_{n-1}}(V^\mu) &\cong& \bigoplus_{\lambda \vdash n,\, \mu \subseteq \lambda,} V^\lambda,
\qquad\hbox{for $\mu \vdash (n-1)$.}
\label{eq:ind}
\end{eqnarray}
In each case $\lambda/\mu$ consists of a single box. 
Starting with the trivial representation $V^{(n)}\cong \CC$ and iterating the restriction (\ref{eq:res}) and induction
(\ref{eq:ind}) rules, we see that  the irreducible $S_n$-representations that appear in
$V^{\otimes k}$ are labeled by the partitions in
\begin{equation}
\Lambda_n^k = \left\{ \lambda \vdash n\, \Big| \, |\lambda| - \lambda_1 \le k\,\right\},
\label{eq:irreducibles}
\end{equation}
and the irreducible $S_{n-1}$-representations that appear in $V^{\otimes k}$ are labeled by the partitions in $\Lambda_{n-1}^k$.

There is an action of $\CC A_k(n)$ on $V^{\otimes k}$ (see \cite{Jo,MR,HR2}) that commutes with
$S_n$ and maps $\CC A_k(n)$ surjectively onto the centralizer $End_{S_n}(V^{\otimes k})$. 
This generalizes the actions of $\CC[S_k]$, by place permutations,
and $B_k(n)$ on $V^{\otimes k}$. 
Furthermore, when $n \ge 2k$ we have
\begin{equation}
\CC A_k(n) \cong End_{S_n}(V^{\otimes k}) \qquad\hbox{and}\qquad
\CC A_{k+{1\over2}}(n) \cong End_{S_{n-1}}(V^{\otimes k}).
\end{equation}
By convention, we let $\CC A_0(n) = \CC A_{\half}(n) = \CC$.
Since anything that commutes with $S_n$ on $V^{\otimes k}$ will also commute
with $S_{n-1}$, we have $\CC A_{k}(n) \subseteq \CC A_{k+{1\over2}}(n)$, and thus
\begin{equation}
\CC A_0(n)  \subseteq \CC A_{{1\over2}}(n) \subseteq \CC A_1(n) \subseteq \CC A_{1{1\over2}}(n) \subseteq \cdots \subseteq \CC A_{(k-{1\over2}}(n) 
\subseteq \CC A_k(n).
\label{eq:tower}
\end{equation}

The {\it Bratteli diagram\/} for $\CC A_k(n)$ consists of rows of vertices, with the rows labeled
by $0, {1\over2}, 1, 1{1\over2}, \ldots, k$, such that the vertices
in row $i$ are $\Lambda_n^i$ and the vertices in row $i+{1\over2}$ are $\Lambda_{n-1}^{i}$.
Two vertices  are connected by an edge if they are in consecutive
rows and they differ by exactly one box. 
Figure \ref{fig:bratteli} shows the Bratteli diagram for $\CC A_k(6)$.  Let $n \ge 2k$.
By  double centralizer theory (see, for example, \cite{HR2}), we know that
\begin{itemize}
\item[(1)] The irreducible representations of $\CC A_k(n)$ can be indexed by $\Lambda_n^k$, so we let 
$M^\lambda_k$ denote the irreducible $\CC A_k(n)$ representation indexed by $\lambda \in \Lambda_n^k$.
\item[(2)]  The decomposition of $V^{\otimes k}$ as an $S_n \times \CC A_k(n)$-bimodule is given by
\begin{equation}
V^{\otimes k} \cong \bigoplus_{\lambda \in \Lambda_n^k}\ V^\lambda \otimes M^\lambda_k.
\label{eq:SchurWeyl}
\end{equation}
\item[(3)]  The dimension of $M^\lambda_k$ equals the multiplicity of $V^\lambda$ in $V^{\otimes k}$.
The  edges in the Bratteli diagram exactly follow the restriction and induction rules in (\ref{eq:res}), and (\ref{eq:ind}). 
and so
$$
m^\lambda_k = \dim(M^\lambda_k) = \left\{ 
\begin{array}{l}
\hbox{ the number of paths from the top }     \\
\hbox{ of the Bratteli diagram to $\lambda$ } \\
\end{array}
\right\}.
$$
\end{itemize}

\begin{figure}
\caption{Bratteli Diagram for $\CC A_k(6)$}
\label{fig:bratteli}
$$ 
{\beginpicture
\setcoordinatesystem units <0.2cm,0.2cm>         % sets scale
\setplotarea x from -6 to 40, y from 4 to 31   % sets plot size up
\linethickness=0.5pt                          % sets line thickness

%%%%%%%%%% level 3 %%%%%%%%%%%%%%%%%%%%%%%%%%%%%%%%%%%%%%%%%%%%
\put{$k = 3:$} at  -6 -.5
\putrectangle corners at 0 -1 and 1 0
\putrectangle corners at 1 -1 and 2 0
\putrectangle corners at 2 -1 and 3 0
\putrectangle corners at 3 -1 and 4 0
\putrectangle corners at 4 -1 and 5 0
\putrectangle corners at 5 -1 and 6 0

\putrectangle corners at 8 -1 and 9 0
\putrectangle corners at 9 -1 and 10 0
\putrectangle corners at 10 -1 and 11 0
\putrectangle corners at 11 -1 and 12 0
\putrectangle corners at 12 -1 and 13 0
\putrectangle corners at 8 -2 and 9 -1

\putrectangle corners at 15 -1 and 16 0
\putrectangle corners at 16 -1 and 17 0
\putrectangle corners at 17 -1 and 18 0
\putrectangle corners at 18 -1 and 19 0
\putrectangle corners at 15 -2 and 16 -1
\putrectangle corners at 16 -2 and 17 -1

\putrectangle corners at 21 -1 and 22 0
\putrectangle corners at 22 -1 and 23 0
\putrectangle corners at 23 -1 and 24 0
\putrectangle corners at 24 -1 and 25 0
\putrectangle corners at 21 -2 and 22 -1
\putrectangle corners at 21 -3 and 22 -2

\putrectangle corners at 27 -1 and 28 0
\putrectangle corners at 28 -1 and 29 0
\putrectangle corners at 29 -1 and 30 0
\putrectangle corners at 27 -2 and 28 -1
\putrectangle corners at 28 -2 and 29 -1
\putrectangle corners at 29 -2 and 30 -1

\putrectangle corners at 32 -1 and 33 0
\putrectangle corners at 33 -1 and 34 0
\putrectangle corners at 34 -1 and 35 0
\putrectangle corners at 32 -2 and 33 -1
\putrectangle corners at 33 -2 and 34 -1
\putrectangle corners at 32 -3 and 33 -2

\putrectangle corners at 37 -1 and 38 0
\putrectangle corners at 38 -1 and 39 0
\putrectangle corners at 39 -1 and 40 0
\putrectangle corners at 37 -2 and 38 -1
\putrectangle corners at 37 -3 and 38 -2
\putrectangle corners at 37 -4 and 38 -3
%-------------------------------------------------

\plot 3 1 3 6  /
\plot 8 1 4 6  /
\plot 9 1 9 5  /
\plot 15 1 10 6  /
\plot 16 1 16 5  /
\plot 27 1 17 5  /
\plot 32 1 18 5  /
\plot 21 1 11 6  /
\plot 23 1 23 5  /
\plot 33 1 24 5  /
\plot 37 1 25 5  /

%%%%%%%%%% level 2+ %%%%%%%%%%%%%%%%%%%%%%%%%%%%%%%%%%%%%%%%%%%%
\put{$k = 2\frac{1}{2}:$} at  -6 7.5
\putrectangle corners at 0 7 and 1 8
\putrectangle corners at 1 7 and 2 8
\putrectangle corners at 2 7 and 3 8
\putrectangle corners at 3 7 and 4 8
\putrectangle corners at 4 7 and 5 8

\putrectangle corners at 8 7  and 9 8
\putrectangle corners at 9 7  and 10 8
\putrectangle corners at 10 7 and 11 8
\putrectangle corners at 11 7 and 12 8
\putrectangle corners at 8  6 and 9 7

\putrectangle corners at 15 7 and 16 8
\putrectangle corners at 16 7 and 17 8
\putrectangle corners at 17 7 and 18 8
\putrectangle corners at 15 6 and 16 7
\putrectangle corners at 16 6 and 17 7

\putrectangle corners at 21 7 and 22 8
\putrectangle corners at 22 7 and 23 8
\putrectangle corners at 23 7 and 24 8
\putrectangle corners at 21 6 and 22 7
\putrectangle corners at 21 5 and 22 6
%-------------------------------------------------

\plot 3 9 3 14  /
\plot 8 13 4 9  /
\plot 9 13 9 9  /
\plot 15 13 10 9  /
\plot 17 13 17 9  /
\plot 20 13 11 9  /
\plot 23 13 23 9  /

%%%%%%%%%% level 2 %%%%%%%%%%%%%%%%%%%%%%%%%%%%%%%%%%%%%%%%%%%%
\put{k = $2:$} at  -6 15.5
\putrectangle corners at 0  15 and 1  16
\putrectangle corners at 1  15 and 2  16
\putrectangle corners at 2  15 and 3  16
\putrectangle corners at 3  15 and 4  16
\putrectangle corners at 4  15 and 5  16
\putrectangle corners at 5  15 and 6  16

\putrectangle corners at 8  15 and 9   16
\putrectangle corners at 9  15 and 10  16
\putrectangle corners at 10 15 and 11  16
\putrectangle corners at 11 15 and 12  16
\putrectangle corners at 12 15 and 13  16
\putrectangle corners at 8  14 and 9  15

\putrectangle corners at 15  15 and 16  16
\putrectangle corners at 16  15 and 17  16
\putrectangle corners at 17  15 and 18  16
\putrectangle corners at 18  15 and 19  16
\putrectangle corners at 15  14 and 16  15
\putrectangle corners at 16  14 and 17  15

\putrectangle corners at 21  15 and 22  16
\putrectangle corners at 22  15 and 23  16
\putrectangle corners at 23  15 and 24  16
\putrectangle corners at 24  15 and 25  16
\putrectangle corners at 21  14 and 22  15
\putrectangle corners at 21  13 and 22 14
%-------------------------------------------------

\plot 3 17 3 22  /
\plot 8 17 4 22  /
\plot 9 17 9 21  /
\plot 15 17 10 21  /
\plot 21 17 11 21  /

%%%%%%%%%% level 1+ %%%%%%%%%%%%%%%%%%%%%%%%%%%%%%%%%%%%%%%%%%%%
\put{$k = 1\frac{1}{2}:$} at  -6 23.5
\putrectangle corners at 0  23 and 1  24
\putrectangle corners at 1  23 and 2  24
\putrectangle corners at 2  23 and 3  24
\putrectangle corners at 3  23 and 4  24
\putrectangle corners at 4  23 and 5  24

\putrectangle corners at 8  23 and 9   24
\putrectangle corners at 9  23 and 10  24
\putrectangle corners at 10 23 and 11  24
\putrectangle corners at 11 23 and 12  24
\putrectangle corners at 8  22 and 9  23
%-------------------------------------------------

\plot 3 25 3 30  /
\plot 8 29 4 25  /
\plot 9 29 9 25  /

%%%%%%%%%% level 1 %%%%%%%%%%%%%%%%%%%%%%%%%%%%%%%%%%%%%%%%%%%%
\put{$k = 1:$} at  -6 31.5
\putrectangle corners at 0  31 and 1  32
\putrectangle corners at 1  31 and 2  32
\putrectangle corners at 2  31 and 3  32
\putrectangle corners at 3  31 and 4  32
\putrectangle corners at 4  31 and 5  32
\putrectangle corners at 5  31 and 6  32

\putrectangle corners at 8  31 and 9   32
\putrectangle corners at 9  31 and 10  32
\putrectangle corners at 10 31 and 11  32
\putrectangle corners at 11 31 and 12  32
\putrectangle corners at 12 31 and 13  32
\putrectangle corners at 8  30 and 9  32

%%%%%%%%%% level 0+ %%%%%%%%%%%%%%%%%%%%%%%%%%%%%%%%%%%%%%%%%%%%
\put{$k = \frac{1}{2}:$} at  -6 39.5
\putrectangle corners at 0  39 and 1  40
\putrectangle corners at 1  39 and 2  40
\putrectangle corners at 2  39 and 3  40
\putrectangle corners at 3  39 and 4  40
\putrectangle corners at 4  39 and 5  40

%-------------------------------------------------

\plot 3 33 3 38  /
\plot 8 33 4 38  /

%%%%%%%%%% level 0 %%%%%%%%%%%%%%%%%%%%%%%%%%%%%%%%%%%%%%%%%%%%
\put{$k = 0:$} at  -6 47.5
\putrectangle corners at 0  47 and 1  48
\putrectangle corners at 1  47 and 2  48
\putrectangle corners at 2  47 and 3  48
\putrectangle corners at 3  47 and 4  48
\putrectangle corners at 4  47 and 5  48
\putrectangle corners at 5  47 and 6  48

%-------------------------------------------------

\plot 3 41 3 46  /

\endpicture}
$$
\end{figure}

\subsection{Vacillating Tableaux}

The dimension of the irreducible $S_n$-module $V^\lambda$ equals the number $f^\lambda$ of
standard tableaux of shape $\lambda$. A standard tableau of shape $\lambda$ is a filling of the Young diagram of $\lambda$ with the numbers $1, 2, \ldots, n$ in such a way that each number appears exactly once, the rows increase from left to right, and the columns increase from top to bottom.  We can identify
a standard tableaux $T_\lambda$ of shape $\lambda$ with a sequence $(\emptyset = \lambda^{(0)}, \lambda^{(1)},
\ldots, \lambda^{(n)} = \lambda)$ such that  $|\lambda^{(i)}| = i$, $\lambda^{(i)} \subseteq \lambda^{(i+1)}$,
and such that $\lambda^{(i)}/\lambda^{(i-1)}$ is the box containing $i$ in $T_\lambda$. For example,
$$
{\beginpicture
\setcoordinatesystem units <0.35cm,0.35cm>         % sets scale
\setplotarea x from 0 to 5, y from -1 to 1   % sets plot size up
\linethickness=0.5pt   
\put{1} at 1.5 .5
\put{2} at 1.5 -.5
\put{3} at 2.5 .5
\put{4} at 3.5 .5
\put{5} at 2.5 -.5
\putrectangle corners at 1  0 and  2 1
\putrectangle corners at 1  -1 and 2 0
\putrectangle corners at 2  0 and  3 1
\putrectangle corners at 3  0 and  4 1
\putrectangle corners at 2  -1 and 3 0
\endpicture}
= \quad
\left(\ \emptyset \ ,
{\beginpicture
\setcoordinatesystem units <0.3cm,0.3cm>         % sets scale
\setplotarea x from .8 to 2.2, y from -1 to 1   % sets plot size up
\linethickness=0.5pt   
\putrectangle corners at 1  0 and  2 1
\endpicture},
{\beginpicture
\setcoordinatesystem units <0.3cm,0.3cm>         % sets scale
\setplotarea x from .8 to 2.2, y from -1 to 1   % sets plot size up
\linethickness=0.5pt   
\putrectangle corners at 1  0 and  2 1
\putrectangle corners at 1  -1 and 2 0
\endpicture},
{\beginpicture
\setcoordinatesystem units <0.3cm,0.3cm>         % sets scale
\setplotarea x from .8 to 2.2, y from -1 to 1   % sets plot size up
\linethickness=0.5pt   
\putrectangle corners at 1  0 and  2 1
\putrectangle corners at 1  -1 and 2 0
\putrectangle corners at 2  0 and  3 1
\endpicture},
{\beginpicture
\setcoordinatesystem units <0.3cm,0.3cm>         % sets scale
\setplotarea x from .8 to 3.2, y from -1 to 1   % sets plot size up
\linethickness=0.5pt   
\putrectangle corners at 1  0 and  2 1
\putrectangle corners at 1  -1 and 2 0
\putrectangle corners at 2  0 and  3 1
\putrectangle corners at 3  0 and  4 1
\endpicture},
{\beginpicture
\setcoordinatesystem units <0.3cm,0.3cm>         % sets scale
\setplotarea x from .8 to 3.2, y from -1 to 1   % sets plot size up
\linethickness=0.5pt   
\putrectangle corners at 1  0 and  2 1
\putrectangle corners at 1  -1 and 2 0
\putrectangle corners at 2  0 and  3 1
\putrectangle corners at 3  0 and  4 1
\putrectangle corners at 2  -1 and 3 0
\endpicture}
\right).
$$
The sequence $(\emptyset = \lambda^{(0)}, \lambda^{(1)}, \ldots, \lambda^{(n)} = \lambda)$ is
a path in Young's lattice, which is the Bratteli diagram for $S_n$.  The number
of standard tableaux $f^\lambda$ can be computed using the hook formula (see \cite[\S3.10]{Sag}).
We let $\SYT(\lambda)$ denote the set of standard tableaux of shape $\lambda$.

Let $\lambda \in \Lambda_n^k$.  A  {\it vacillating tableaux of shape $\lambda$ and length $2k$} is a
sequence of partitions,
$$
\left((n) = \lambda^{(0)}, \lambda^{({1\over2})}, \lambda^{(1)}, \lambda^{(1{1\over2})}, \ldots, 
\lambda^{(k-{1\over2})},\lambda^{(k)} = \lambda\right),
$$
satisfying, for each $i$,
\begin{itemize}
\item[(1)] $\lambda^{(i)} \in \Lambda_n^i$ and  $\lambda^{(i+{1\over2})} \in \Lambda_{n-1}^{i},$
\item[(2)] $\lambda^{(i)} \supseteq  \lambda^{(i+{1\over2})}$ and $|\lambda^{(i)} / \lambda^{(i+{1\over2})}| = 1,$
\item[(3)] $\lambda^{(i+{1\over2})} \subseteq \lambda^{(i+1)}$ and $|\lambda^{(i+1)} / \lambda^{(i+{1\over2})}| = 1.$
\end{itemize}
The vacillating tableaux of shape $\lambda$ correspond exactly with paths from the top of the Bratteli diagram to $\lambda$. Thus, if we let $\PT_k(\lambda)$ denote the set
of vacillating tableaux of shape $\lambda$ and length $k$, then 
\begin{equation}
m^\lambda_k = \dim(M^\lambda_k) = | \PT_k(\lambda) |.
\end{equation}

Let $n \ge 2k$, and for a partition $\lambda$, define $\lambda^\ast$ and
$\bar
\lambda$ as follows,
\[
\begin{array}{llll}
\hbox{if} & \lambda = (\lambda_1, \ldots, \lambda_\ell)\ \vdash \
n, &
\hbox{then} & \lambda^\ast = (\lambda_2, \ldots, \lambda_\ell)\ \vdash\
(n - \lambda_1), \\
\hbox{if} & \lambda = (\lambda_1, \ldots, \lambda_\ell)\ \vdash\
s \le n, &
\hbox{then} & \bar\lambda = (n - s, \lambda_1, \ldots, \lambda_\ell)
\ \vdash \ n. \\
\end{array}
\]
Since $n \ge 2k$ we are  guaranteed that $\bar\lambda$ is a partition
and that $0 \le |\lambda^\ast| \le k$.  The sets
\begin{equation}
\Lambda_n^k = \left\{ \, \lambda \vdash n\, \Big| \,  |\lambda^\ast| \le k\,\right\}
\qquad\hbox{and}\qquad
\Gamma_k = \left\{ \, \lambda \vdash t\, \Big| \, 0 \le t \le k\, \right\}
\end{equation}
are in bijection with one another using the  maps,
\begin{equation}
\begin{array} {ccc}
\Lambda_{n}^k &\to& \Gamma_k \\
\lambda & \mapsto & \lambda^\ast \\
\end{array}
\qquad\hbox{and}\qquad
\begin{array} {ccc}
\Gamma_k &\to& \Lambda_{n}^k \\
\lambda & \mapsto & \bar \lambda \\
\end{array}. 
\label{eq:bijection}
\end{equation}
Via these bijections, we can use either $\Gamma_k$ or $\Lambda_n^k$ to
index the irreducible representations of $\CC A_k(n)$.
 For example, the following sequences represent
the same vacillating tableau $P_\lambda$, the first using diagrams from $\Lambda_n^k$
and the second from $\Gamma_k$,
\begin{eqnarray*}
P_\lambda &=&
\left(
{\beginpicture
\setcoordinatesystem units <0.25cm,0.25cm> % sets scale
\setplotarea x from 0 to 0, y from -1 to 1    % sets plot size up
\linethickness=0.5pt      
\putrectangle corners at 0 1 and 1 0
\putrectangle corners at 1 1 and 2 0
\putrectangle corners at 2 1 and 3 0
\putrectangle corners at 3 1 and 4 0
\putrectangle corners at 4 1 and 5 0
\putrectangle corners at 5 1 and 6 0
\endpicture},
{\beginpicture
\setcoordinatesystem units <0.25cm,0.25cm> % sets scale
\setplotarea x from 0 to 0, y from -1 to 1    % sets plot size up
\linethickness=0.5pt      
\putrectangle corners at 0 1 and 1 0
\putrectangle corners at 1 1 and 2 0
\putrectangle corners at 2 1 and 3 0
\putrectangle corners at 3 1 and 4 0
\putrectangle corners at 4 1 and 5 0
\endpicture},
{\beginpicture
\setcoordinatesystem units <0.25cm,0.25cm> % sets scale
\setplotarea x from 0 to 0, y from -1 to 1    % sets plot size up
\linethickness=0.5pt      
\putrectangle corners at 0 1 and 1 0
\putrectangle corners at 1 1 and 2 0
\putrectangle corners at 2 1 and 3 0
\putrectangle corners at 3 1 and 4 0
\putrectangle corners at 4 1 and 5 0
\putrectangle corners at 0 0 and 1 -1
\endpicture},
{\beginpicture
\setcoordinatesystem units <0.25cm,0.25cm> % sets scale
\setplotarea x from 0 to 0, y from -1 to 1    % sets plot size up
\linethickness=0.5pt      
\putrectangle corners at 0 1 and 1 0
\putrectangle corners at 1 1 and 2 0
\putrectangle corners at 2 1 and 3 0
\putrectangle corners at 3 1 and 4 0
\putrectangle corners at 0 0 and 1 -1
\endpicture},
{\beginpicture
\setcoordinatesystem units <0.25cm,0.25cm> % sets scale
\setplotarea x from 0 to 0, y from -1 to 1    % sets plot size up
\linethickness=0.5pt      
\putrectangle corners at 0 1 and 1 0
\putrectangle corners at 1 1 and 2 0
\putrectangle corners at 2 1 and 3 0
\putrectangle corners at 3 1 and 4 0
\putrectangle corners at 0 0 and 1 -1
\putrectangle corners at 1 0 and 2 -1
\endpicture},
{\beginpicture
\setcoordinatesystem units <0.25cm,0.25cm> % sets scale
\setplotarea x from 0 to 0, y from -1 to 1    % sets plot size up
\linethickness=0.5pt      
\putrectangle corners at 0 1 and 1 0
\putrectangle corners at 1 1 and 2 0
\putrectangle corners at 2 1 and 3 0
\putrectangle corners at 1 0 and 2 -1
\putrectangle corners at 0 0 and 1 -1
\endpicture},
{\beginpicture
\setcoordinatesystem units <0.25cm,0.25cm> % sets scale
\setplotarea x from 0 to 0, y from -1 to 1    % sets plot size up
\linethickness=0.5pt      
\putrectangle corners at 0 1 and 1 0
\putrectangle corners at 1 1 and 2 0
\putrectangle corners at 2 1 and 3 0
\putrectangle corners at 3 1 and 4 0
\putrectangle corners at 0 0 and 1 -1
\putrectangle corners at 1 0 and 2 -1
\endpicture}
\right), \\
&=& \left(
\emptyset,
\emptyset,
{\beginpicture
\setcoordinatesystem units <0.25cm,0.25cm> % sets scale
\setplotarea x from 0 to 0, y from -1 to 1    % sets plot size up
\linethickness=0.5pt      
\putrectangle corners at 0 1 and 1 0
\endpicture},
{\beginpicture
\setcoordinatesystem units <0.25cm,0.25cm> % sets scale
\setplotarea x from 0 to 0, y from -1 to 1    % sets plot size up
\linethickness=0.5pt      
\putrectangle corners at 0 1 and 1 0
\endpicture},
{\beginpicture
\setcoordinatesystem units <0.25cm,0.25cm> % sets scale
\setplotarea x from 0 to 0, y from -1 to 1    % sets plot size up
\linethickness=0.5pt      
\putrectangle corners at 0 1 and 1 0
\putrectangle corners at 1 1 and 2 0
\endpicture},
{\beginpicture
\setcoordinatesystem units <0.25cm,0.25cm> % sets scale
\setplotarea x from 0 to 0, y from -1 to 1    % sets plot size up
\linethickness=0.5pt      
\putrectangle corners at 0 1 and 1 0
\putrectangle corners at 1 1 and 2 0
\endpicture},
{\beginpicture
\setcoordinatesystem units <0.25cm,0.25cm> % sets scale
\setplotarea x from 0 to 0, y from -1 to 1    % sets plot size up
\linethickness=0.5pt      
\putrectangle corners at 0 1 and 1 0
\putrectangle corners at 1 1 and 2 0
\endpicture},
\right).
\end{eqnarray*}
For our bijection in Section 3 we will use $\Lambda_n^k$, and for our
bijection in Section 4 we will use $\Gamma_k$.    The Bratteli diagram for $\CC A_k(n), n \ge 2k$, is  given 
in Figure 2, using labels from $\Gamma_k$, along with the number of vacillating tableaux for each shape $\lambda$. 

\begin{figure}
\label{fig:BratteliPartitionGeneric}
\caption{Bratteli Diagram for $\CC A_k(n), n \ge 2k$.}
$$
{\beginpicture
\setcoordinatesystem units <0.175cm,0.175cm>         % sets scale
\setplotarea x from -8 to 40, y from 4 to 31   % sets plot size up
\linethickness=0.5pt                          % sets line thickness
%%%%%%%%%% level 3 %%%%%%%%%%%%%%%%%%%%%%%%%%%%%%%%%%%%%%%%%%%%
\put{$k = 3:$} at  -6 -.5
\put{$\emptyset$} at 3 -.5
\putrectangle corners at 8.5 -1 and 9.5 0
\putrectangle corners at 15 -1 and 16 0
\putrectangle corners at 16 -1 and 17 0
\putrectangle corners at 22 -1 and 23 0
\putrectangle corners at 22 -2 and 23 -1
\putrectangle corners at 27 -1 and 28 0
\putrectangle corners at 28 -1 and 29 0
\putrectangle corners at 26 -1 and 27 0
\putrectangle corners at 32 -1 and 33 0
\putrectangle corners at 33 -1 and 34 0
\putrectangle corners at 32 -2 and 33 -1
\putrectangle corners at 37 -1 and 38 0
\putrectangle corners at 37 -2 and 38 -1
\putrectangle corners at 37 -3 and 38 -2
%-------------------------------------------------
%
\plot 3 1 3 6  /
\plot 8 1 4 6  /
\plot 9 1 9 6  /
\plot 15 1 10 6  /
\plot 16 1 16 6  /
\plot 26 1 17 6  /
\plot 31.5 1 18 6  /
\plot 21 1 11 6  /
\plot 22.5 1 22.5 5.5  /
\plot 32.5 1 23 5.5  /
\plot 37 1 24 5.5  /
%
%%%%%%%%%% level 2+ %%%%%%%%%%%%%%%%%%%%%%%%%%%%%%%%%%%%%%%%%%%%
\put{$k = 2{1\over2}:$} at  -6 7.5
\put{$\emptyset$} at 3 7.5
\putrectangle corners at 8.5 7  and 9.5 8
\putrectangle corners at 15 7 and 16 8
\putrectangle corners at 16 7 and 17 8
\putrectangle corners at 22 7 and 23 8
\putrectangle corners at 22 6 and 23 7
%-------------------------------------------------
%
\plot 3 9 3 13  /
\plot 8 13 4 9  /
\plot 9 13 9 9  /
\plot 15 13 10 9  /
\plot 16 13 16 9  /
\plot 22 12.5 11 9  /
\plot 22.5 12.5 22.5 9  /
%
%
%%%%%%%%%% level 2 %%%%%%%%%%%%%%%%%%%%%%%%%%%%%%%%%%%%%%%%%%%%
\put{k = $2:$} at  -6 14.5
\put{$\emptyset$} at 3 14.5
\putrectangle corners at 8  14 and 9   15
\putrectangle corners at 15  14 and 16  15
\putrectangle corners at 16  14 and 17  15
\putrectangle corners at 22  14 and 23  15
\putrectangle corners at 22  13 and 23  14
%-------------------------------------------------
%
\plot 3 16 3 20  /
\plot 8 16 4 20  /
\plot 9 16 9 20  /
\plot 15 16 9.5 20  /
\plot 22 16 10 20  /
%
%
%%%%%%%%%% level 1+ %%%%%%%%%%%%%%%%%%%%%%%%%%%%%%%%%%%%%%%%%%%%
\put{$k = 1{1\over2}:$} at  -6 21.5
\put{$\emptyset$} at 3 21.5
\putrectangle corners at 8  21 and 9   22
%-------------------------------------------------
%
\plot 3 23 3 27  /
\plot 8 27 4 23  /
\plot 9 27 9 23  /
%
%
%%%%%%%%%% level 1 %%%%%%%%%%%%%%%%%%%%%%%%%%%%%%%%%%%%%%%%%%%%
\put{$k = 1:$} at  -6 28.5
\put{$\emptyset$} at  3 28.5
\putrectangle corners at 8  28 and 9  29
%
%%%%%%%%%% level 0+ %%%%%%%%%%%%%%%%%%%%%%%%%%%%%%%%%%%%%%%%%%%%
\put{$k = {1\over2}:$} at  -6 35.5
\put{$\emptyset$} at  3 35.5
%
%-------------------------------------------------
%
\plot 3 30 3 34  /
\plot 8 30 4 34  /
%
%
%%%%%%%%%% level 0 %%%%%%%%%%%%%%%%%%%%%%%%%%%%%%%%%%%%%%%%%%%%
\put{$k = 0:$} at  -6 42.5
\put{$\emptyset$} at  3 42.5
%
%-------------------------------------------------
%
\plot 3 37 3 41  /
\endpicture} 
{\beginpicture
\setcoordinatesystem units <0.175cm,0.175cm>         % sets scale
\setplotarea x from 0 to 40, y from 4 to 31   % sets plot size up
\linethickness=0.5pt                          % sets line thickness
%%%%%%%%%% level 3 %%%%%%%%%%%%%%%%%%%%%%%%%%%%%%%%%%%%%%%%%%%%
\put{5} at 3 -.5
\put{10}  at 8.5 -.5 
\put{6} at 15.5 -.5 
\put{6} at 22.5 -.5 
\put{1} at 27.5 -.5 
\put{2} at 32.5 -.5 
\put{1} at 37.5 -.5 
\plot 3 1 3 6  /
\plot 8 1 4 6  /
\plot 9 1 9 6  /
\plot 15 1 10 6  /
\plot 15.5 1 15.5 6  /
\plot 26 1 17 6  /
\plot 31.5 1 18 6  /
\plot 21 1 11 6  /
\plot 22.5 1 22.5 6  /
\plot 32.5 1 23 6  /
\plot 37 1 24 6  /
%
%%%%%%%%%% level 2+ %%%%%%%%%%%%%%%%%%%%%%%%%%%%%%%%%%%%%%%%%%%%
\put{5} at 3 7.5
\put{5}  at 8.5 7.5  
\put{1} at 15.5 7.5 
\put{1} at 22.5 7.5 
\plot 3 9 3 13  /
\plot 8 13 4 9  /
\plot 9 13 9 9  /
\plot 15 13 10 9  /
\plot 15.5 13 15.5 9  /
\plot 22 13 11 9  /
\plot 22.5 13 22.5 9  /
% level 2 %%%%%%%%%%%%%%%%%%%%
\put{2} at 3 14.5
\put{3} at 8.5  14.5 
\put{1} at 15.5  14.5 
\put{1} at 22.5  14.5 
\plot 3 16 3 20  /
\plot 8 16 4 20  /
\plot 9 16 9 20  /
\plot 15 16 9.5 20  /
\plot 22 16 10 20  /
%level 1+ %%%%%%%%%%%%%%%%%%
\put{2} at 3 21.5
\put{1} at 8.5  21.5 
\plot 3 23 3 27  /
\plot 8 27 4 23  /
\plot 9 27 9 23  /
% level 1 %%%%%%%%%%%%%%%%%%\
\put{1} at  3 28.5
\put{1} at 8.5  28.5 
%level 0+ %%%%%%%%%%%%%%
\put{1} at  3 35.5
\plot 3 30 3 34  /
\plot 8 30 4 34  /
% level 0 %%%%%%%%%%%%%%
\put{1} at  3 42.5
\plot 3 37 3 41  /
\endpicture} 
$$
\end{figure}

\section{A Bijective Proof of $n^k = \displaystyle{\sum_{\lambda \in
\Lambda_n^k}}\, f^\lambda\, m^\lambda_k$}

Comparing dimensions on both sides of the identity  (\ref{eq:SchurWeyl}) gives
\begin{equation}
n^k = \sum_{\lambda \in \Lambda_n^k}\, f^\lambda\, m^\lambda_k, \qquad
\hbox{ for $n \ge 2k$.}
\label{eq:identity1}
\end{equation}
For example, when $n = 6$ and $k = 3$, the  $f^\lambda$ in the bottom row of Figure \ref{fig:bratteli} (see also Figure 2) are 1, 5, 9, 10, 5, 16, 10, the corresponding $m^\lambda_3$ are 5, 10, 6, 6, 1, 2, 1, and we have
$ 6^3 = 216 = 1\!\cdot\!5 + 5 \!\cdot\! 10 + 9 \!\cdot\! 6 + 10 \!\cdot\! 6 + 
5 \!\cdot\! 1 + 16 \!\cdot\! 2 + 10\!\cdot\! 1.$

To give a combinatorial proof of (\ref{eq:identity1}), we need to find a bijection
of the form
\begin{equation}
\left\{\, (i_1, \ldots, i_k) \, \Big|\, 1 \le i_j \le n\,\right\}
\longleftrightarrow
\bigsqcup_{\lambda \in \Lambda_n^k} \SYT(\lambda) \times \PT_k(\lambda).
\label{eq:identity1A}
\end{equation}
To do so, we construct an invertible function that turns a sequence $(i_1, \ldots, i_k)$ of numbers
in the range $1 \le i_j \le n$ into a pair $(T_\lambda, P_\lambda)$ consisting of a standard  tableaux $T_\lambda$
of shape $\lambda$ and vacillating tableaux $P_\lambda$ of shape $\lambda$ and length $2k$ for some $\lambda \in
\Lambda_n^k$.  Our bijection uses  jeu de taquin and RSK column insertion.

If $T$ is a standard tableau of shape $\lambda \vdash n$, then Sch\" utzenberger's \cite{Scu}
{\it jeu de taquin} provides an algorithm for removing the box containing $x$ from
$T$ and producing a standard tableau $S$ of shape $\mu \vdash (n-1)$ with $\mu \subseteq \lambda$ and entries $\{1,
\ldots, n\}\setminus\{x\}.$  We only need a special case of jeu de taquin for our purposes. See
\cite[\S 3.7]{Sag} or \cite[\S A1.2]{Sta} for the full-strength version and its applications.

If $S$ is a standard tableau, let $S_{i,j}$ denote the entry of $S$ in row $i$ (numbered left-to-right) and column
$j$ (numbered top-to-bottom).  We say that a {\it corner\/} of $S$ is a box whose removal leaves the Young diagram of
a partition. Thus the corners of $S$ are the boxes that are both at the end of a row and the end of a column.
The following algorithm will delete $x$ from $T$ leaving a standard tableau $S$ with $x$ removed.
We denote this process by $x \jdt T$.
\begin{tabbing}
\quad \= \qquad\=\quad \kill
\> 1. Let $c = S_{i,j}$ be the box containing $x$. \\
\> 2. \keyw{While} $c$ is not a corner,  \keyw{do} \\
\>\> A. Let $c'$ be the box containing $\min\{S_{i+1,j},S_{i,j+1}\}$; \\
\>\>B. Exchange the positions of $c$ and $c'$.\\
\> 3. Delete $c$.
\end{tabbing}
If only one of $S_{i+1,j},S_{i,j+1}$ exists at step 2.A, then the minimum is
taken to be that single value.  Below is an example of  $2 \jdt T$, 
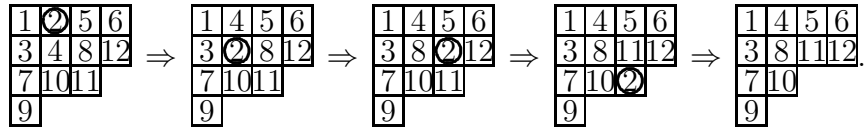
$$
{\beginpicture
\setcoordinatesystem units <0.4cm,0.4cm>         % sets scale
\setplotarea x from 1 to 5, y from -1 to 1   % sets plot size up
\linethickness=0.5pt   
\putrectangle corners at 1  0 and  2 1
\putrectangle corners at 2  0 and  3 1
\putrectangle corners at 3  0 and  4 1
\putrectangle corners at 4  0 and  5 1
\putrectangle corners at 1  1 and  2 2
\putrectangle corners at 2  1 and  3 2
\putrectangle corners at 3  1 and  4 2
\putrectangle corners at 4  1 and  5 2
\putrectangle corners at 1  -1 and  2 0
\putrectangle corners at 2  -1 and  3 0
\putrectangle corners at 3  -1 and  4 0
\putrectangle corners at 1  -2 and  2 -1
\put{1} at 1.5 1.5
\put{2} at 2.5 1.5
\circulararc 360 degrees from 2.1 1.5 center at 2.5 1.5
\put{3} at 1.5 0.5
\put{5} at 3.5 1.5
\put{4} at 2.5 .5
\put{6} at 4.5 1.5
\put{7} at 1.5 -0.5
\put{10} at 2.5 -0.5
\put{8} at 3.5 0.5
\put{9} at 1.5 -1.5
\put{12} at 4.5 .5
\put{11} at 3.5 -0.5
\endpicture} \,\Rightarrow\,
{\beginpicture
\setcoordinatesystem units <0.4cm,0.4cm>         % sets scale
\setplotarea x from 1 to 5, y from -1 to 1   % sets plot size up
\linethickness=0.5pt   
\putrectangle corners at 1  0 and  2 1
\putrectangle corners at 2  0 and  3 1
\putrectangle corners at 3  0 and  4 1
\putrectangle corners at 4  0 and  5 1
\putrectangle corners at 1  1 and  2 2
\putrectangle corners at 2  1 and  3 2
\putrectangle corners at 3  1 and  4 2
\putrectangle corners at 4  1 and  5 2
\putrectangle corners at 1  -1 and  2 0
\putrectangle corners at 2  -1 and  3 0
\putrectangle corners at 3  -1 and  4 0
\putrectangle corners at 1  -2 and  2 -1
\put{1} at 1.5 1.5
\put{4} at 2.5 1.5
\put{3} at 1.5 0.5
\put{5} at 3.5 1.5
\put{2} at 2.5 .5
\circulararc 360 degrees from 2.1 0.6 center at 2.5 0.5
\put{6} at 4.5 1.5
\put{7} at 1.5 -0.5
\put{10} at 2.5 -0.5
\put{8} at 3.5 0.5
\put{9} at 1.5 -1.5
\put{12} at 4.5 .5
\put{11} at 3.5 -0.5
\endpicture} \,\Rightarrow\,
{\beginpicture
\setcoordinatesystem units <0.4cm,0.4cm>         % sets scale
\setplotarea x from 1 to 5, y from -1 to 1   % sets plot size up
\linethickness=0.5pt   
\putrectangle corners at 1  0 and  2 1
\putrectangle corners at 2  0 and  3 1
\putrectangle corners at 3  0 and  4 1
\putrectangle corners at 4  0 and  5 1
\putrectangle corners at 1  1 and  2 2
\putrectangle corners at 2  1 and  3 2
\putrectangle corners at 3  1 and  4 2
\putrectangle corners at 4  1 and  5 2
\putrectangle corners at 1  -1 and  2 0
\putrectangle corners at 2  -1 and  3 0
\putrectangle corners at 3  -1 and  4 0
\putrectangle corners at 1  -2 and  2 -1
\put{1} at 1.5 1.5
\put{4} at 2.5 1.5
\put{3} at 1.5 0.5
\put{5} at 3.5 1.5
\put{8} at 2.5 .5
\put{6} at 4.5 1.5
\put{7} at 1.5 -0.5
\put{10} at 2.5 -0.5
\put{2} at 3.5 0.5
\circulararc 360 degrees from 3.1 0.6 center at 3.5 0.5
\put{9} at 1.5 -1.5
\put{12} at 4.5 .5
\put{11} at 3.5 -0.5
\endpicture}  
\,\Rightarrow\,
{\beginpicture
\setcoordinatesystem units <0.4cm,0.4cm>         % sets scale
\setplotarea x from 1 to 5, y from -1 to 1   % sets plot size up
\linethickness=0.5pt   
\putrectangle corners at 1  0 and  2 1
\putrectangle corners at 2  0 and  3 1
\putrectangle corners at 3  0 and  4 1
\putrectangle corners at 4  0 and  5 1
\putrectangle corners at 1  1 and  2 2
\putrectangle corners at 2  1 and  3 2
\putrectangle corners at 3  1 and  4 2
\putrectangle corners at 4  1 and  5 2
\putrectangle corners at 1  -1 and  2 0
\putrectangle corners at 2  -1 and  3 0
\putrectangle corners at 3  -1 and  4 0
\putrectangle corners at 1  -2 and  2 -1
\put{1} at 1.5 1.5
\put{4} at 2.5 1.5
\put{3} at 1.5 0.5
\put{5} at 3.5 1.5
\put{8} at 2.5 .5
\put{6} at 4.5 1.5
\put{7} at 1.5 -0.5
\put{10} at 2.5 -0.5
\put{11} at 3.5 0.5
\put{9} at 1.5 -1.5
\put{12} at 4.5 .5
\put{2} at 3.5 -0.5
\circulararc 360 degrees from 3.9 -0.4 center at 3.5 -0.5
\endpicture} \,\Rightarrow\,
{\beginpicture
\setcoordinatesystem units <0.4cm,0.4cm>         % sets scale
\setplotarea x from 1 to 5, y from -1 to 1   % sets plot size up
\linethickness=0.5pt   
\putrectangle corners at 1  0 and  2 1
\putrectangle corners at 2  0 and  3 1
\putrectangle corners at 3  0 and  4 1
\putrectangle corners at 4  0 and  5 1
\putrectangle corners at 1  1 and  2 2
\putrectangle corners at 2  1 and  3 2
\putrectangle corners at 3  1 and  4 2
\putrectangle corners at 4  1 and  5 2
\putrectangle corners at 1  -1 and  2 0
\putrectangle corners at 2  -1 and  3 0
\putrectangle corners at 1  -2 and  2 -1
\put{1} at 1.5 1.5
\put{4} at 2.5 1.5
\put{3} at 1.5 0.5
\put{5} at 3.5 1.5
\put{8} at 2.5 .5
\put{6} at 4.5 1.5
\put{7} at 1.5 -0.5
\put{10} at 2.5 -0.5
\put{11} at 3.5 0.5
\put{9} at 1.5 -1.5
\put{12} at 4.5 .5
\endpicture}.
$$

A bijective proof of the $S_n$ identity
$
k! = \sum_{\lambda \vdash k} (f^\lambda)^2
$
was originally found by Robinson \cite{Rob} and later found, independently and in the form
we present here, by Schensted \cite{Sch}.  Knuth \cite{Kn} analyzed this algorithm 
and extended it to prove the identity (\ref{eq:Skids}.a). 
See \cite[\S7 Notes]{Sta} for a nice history of the RSK algorithm.
Let $S$ be a tableau of partition
shape
$\mu$,  with $|\mu| < n$, with increasing rows and columns, and with distinct entries from $\{1, \ldots, n\}$. 
Let $x$ be a positive integer that is not in $S$. The following algorithm inserts $x$ into $S$
producing a standard tableau $T$ of shape $\lambda$ with $\mu \subseteq \lambda$, $|\lambda/\mu| = 1$,  whose entries are the union of those from $S$ and $\{x\}$.
We denote this process by $x \RS S$.
\begin{tabbing}
\quad \= \qquad\=\quad \kill
\> 1. Let $R$ be the first row of $S$. \\
\> 2. \keyw{While} $x$ is less than some element in $R$, \keyw{do}  \\
\>\> A. Let $y$ be the smallest element of $R$ greater than $x$; \\
\>\> B. Replace $y \in R$ with $x$; \\
\>\> C. Let $x :=  y$ and let $R$ be the next row.\\
\> 3.  Place $x$ at the end of $R$ (which is possibly empty).
\end{tabbing}
For example, here is the insertion of $2$ into the output of the jeu de taquin example above$$
{\beginpicture
\setcoordinatesystem units <0.4cm,0.4cm>         % sets scale
\setplotarea x from 1 to 5, y from -1 to 1   % sets plot size up
\linethickness=0.5pt   
\putrectangle corners at 1  0 and  2 1
\putrectangle corners at 2  0 and  3 1
\putrectangle corners at 3  0 and  4 1
\putrectangle corners at 4  0 and  5 1
\putrectangle corners at 1  1 and  2 2
\putrectangle corners at 2  1 and  3 2
\putrectangle corners at 3  1 and  4 2
\putrectangle corners at 4  1 and  5 2
\putrectangle corners at 1  -1 and  2 0
\putrectangle corners at 2  -1 and  3 0
\putrectangle corners at 1  -2 and  2 -1
\put{1} at 1.5 1.5
\put{4} at 2.5 1.5
\put{3} at 1.5 0.5
\put{5} at 3.5 1.5
\put{8} at 2.5 .5
\put{6} at 4.5 1.5
\put{7} at 1.5 -0.5
\put{10} at 2.5 -0.5
\put{11} at 3.5 0.5
\put{9} at 1.5 -1.5
\put{12} at 4.5 .5
\put{$\leftarrow$} at 5.6 1.5
\put{2} at 6.5 1.5
\endpicture} \,\Rightarrow\,
{\beginpicture
\setcoordinatesystem units <0.4cm,0.4cm>         % sets scale
\setplotarea x from 1 to 5, y from -1 to 1   % sets plot size up
\linethickness=0.5pt   
\putrectangle corners at 1  0 and  2 1
\putrectangle corners at 2  0 and  3 1
\putrectangle corners at 3  0 and  4 1
\putrectangle corners at 4  0 and  5 1
\putrectangle corners at 1  1 and  2 2
\putrectangle corners at 2  1 and  3 2
\putrectangle corners at 3  1 and  4 2
\putrectangle corners at 4  1 and  5 2
\putrectangle corners at 1  -1 and  2 0
\putrectangle corners at 2  -1 and  3 0
\putrectangle corners at 1  -2 and  2 -1
\put{1} at 1.5 1.5
\put{2} at 2.5 1.5
\put{3} at 1.5 0.5
\put{5} at 3.5 1.5
\put{8} at 2.5 .5
\put{6} at 4.5 1.5
\put{7} at 1.5 -0.5
\put{10} at 2.5 -0.5
\put{11} at 3.5 0.5
\put{9} at 1.5 -1.5
\put{12} at 4.5 .5
\put{$\leftarrow$} at 5.6 0.5
\put{4} at 6.5 0.5
\endpicture}
\,\Rightarrow\,
{\beginpicture
\setcoordinatesystem units <0.4cm,0.4cm>         % sets scale
\setplotarea x from 1 to 5, y from -1 to 1   % sets plot size up
\linethickness=0.5pt   
\putrectangle corners at 1  0 and  2 1
\putrectangle corners at 2  0 and  3 1
\putrectangle corners at 3  0 and  4 1
\putrectangle corners at 4  0 and  5 1
\putrectangle corners at 1  1 and  2 2
\putrectangle corners at 2  1 and  3 2
\putrectangle corners at 3  1 and  4 2
\putrectangle corners at 4  1 and  5 2
\putrectangle corners at 1  -1 and  2 0
\putrectangle corners at 2  -1 and  3 0
\putrectangle corners at 1  -2 and  2 -1
\put{1} at 1.5 1.5
\put{2} at 2.5 1.5
\put{3} at 1.5 0.5
\put{5} at 3.5 1.5
\put{4} at 2.5 .5
\put{6} at 4.5 1.5
\put{7} at 1.5 -0.5
\put{10} at 2.5 -0.5
\put{11} at 3.5 0.5
\put{9} at 1.5 -1.5
\put{12} at 4.5 .5
\put{$\leftarrow$} at 3.6 -0.5
\put{8} at 4.5 -0.5
\endpicture}
\,\Rightarrow\,
{\beginpicture
\setcoordinatesystem units <0.4cm,0.4cm>         % sets scale
\setplotarea x from 1 to 5, y from -1 to 1   % sets plot size up
\linethickness=0.5pt   
\putrectangle corners at 1  0 and  2 1
\putrectangle corners at 2  0 and  3 1
\putrectangle corners at 3  0 and  4 1
\putrectangle corners at 4  0 and  5 1
\putrectangle corners at 1  1 and  2 2
\putrectangle corners at 2  1 and  3 2
\putrectangle corners at 3  1 and  4 2
\putrectangle corners at 4  1 and  5 2
\putrectangle corners at 1  -1 and  2 0
\putrectangle corners at 2  -1 and  3 0
\putrectangle corners at 1  -2 and  2 -1
\put{1} at 1.5 1.5
\put{2} at 2.5 1.5
\put{3} at 1.5 0.5
\put{5} at 3.5 1.5
\put{4} at 2.5 .5
\put{6} at 4.5 1.5
\put{7} at 1.5 -0.5
\put{8} at 2.5 -0.5
\put{11} at 3.5 0.5
\put{9} at 1.5 -1.5
\put{12} at 4.5 .5
\put{$\leftarrow$} at 2.6 -1.5
\put{10} at 3.5 -1.5
\endpicture}
\,\Rightarrow\,
{\beginpicture
\setcoordinatesystem units <0.4cm,0.4cm>         % sets scale
\setplotarea x from 1 to 5, y from -1 to 1   % sets plot size up
\linethickness=0.5pt   
\putrectangle corners at 1  0 and  2 1
\putrectangle corners at 2  0 and  3 1
\putrectangle corners at 3  0 and  4 1
\putrectangle corners at 4  0 and  5 1
\putrectangle corners at 1  1 and  2 2
\putrectangle corners at 2  1 and  3 2
\putrectangle corners at 3  1 and  4 2
\putrectangle corners at 4  1 and  5 2
\putrectangle corners at 1  -1 and  2 0
\putrectangle corners at 2  -1 and  3 0
\putrectangle corners at 1  -2 and  2 -1
\putrectangle corners at 2  -2 and  3 -1
\put{1} at 1.5 1.5
\put{2} at 2.5 1.5
\put{3} at 1.5 0.5
\put{5} at 3.5 1.5
\put{4} at 2.5 .5
\put{7} at 1.5 -0.5
\put{8} at 2.5 -0.5
\put{11} at 3.5 0.5
\put{9} at 1.5 -1.5
\put{10} at 2.5 -1.5
\put{12} at 4.5 .5
\put{6} at 4.5 1.5
\endpicture}
$$
It is possible to invert the process of row insertion using row {\it uninsertion}. See \cite{Sag}
for details.  In the example above, the number 2 and the leftmost tableau are the result of uninserting 10 from the
rightmost tableau.

Given $i_1, \ldots, i_k$, with $1 \le i_j \le n$, we will produce a pair $(T_\lambda, P_\lambda)$, $\lambda \in \Lambda_n^k$,
consisting of a standard tableau $T_\lambda$ and a vacillating tableau $P_\lambda$.  First, initialize the 0th tableau to be the standard tableaux of shape $(n)$, namely,
\begin{equation}
T^{(0)} {\beginpicture
\setcoordinatesystem units <0.5cm,0.5cm>         % sets scale
\setplotarea x from 1 to 7.3, y from -1 to 1   % sets plot size up
\linethickness=0.5pt   
\putrectangle corners at 1  -.5 and  2 .5
\putrectangle corners at 2  -.5 and  3 .5
\putrectangle corners at 3  -.5 and  6 .5
\putrectangle corners at 6  -.5 and  7 .5
\put{1} at 1.5 0
\put{2} at 2.5 0
\put{$\cdots$} at 4.5 0
\put{$n$} at 6.5 0
\endpicture}.
\end{equation}
Then, recursively define standard tableaux $T^{(j+{1\over2})}$ and $T^{(j+1)}$ by
\begin{equation}
\begin{array}{rcl}
T^{(j+{1\over2})} &=& \left( i_{j+1} \jdt T^{(j)} \right), \\
T^{(j+1)} &=& \left( i_{j+1} \RS T^{(j+\frac{1}{2})} \right),  \\
\end{array} \qquad   0 \le j \le k-1.
\end{equation}
Let $\lambda^{(j)} \in \Lambda_n^j$ be the shape of $T^{(j)}$, and let
$\lambda^{(j+{1\over2})} \in \Lambda_{n-1}^j$ be the shape of $T^{(j+{1\over2})}$.  
Then let
\begin{equation}
P_\lambda = \left(\lambda^{(0)}, \lambda^{(\half)}, \lambda^{(1)}, \lambda^{(1\half)},  \ldots, \lambda^{(k)}\right)
\qquad
\hbox{and}
\qquad
T_\lambda = T^{(k)},
\end{equation}
so that $P_\lambda$ is a vacillating tableau of shape $\lambda= \lambda^{(k)} \in \Lambda_n^k$, and
$T_\lambda$ is a standard  tableaux of the same shape $\lambda$.  We denote this iterative ``delete-insert"
process  that associates the pair $(T_\lambda, P_\lambda)$ to the sequence $(i_1, \ldots, i_k)$ by
\begin{equation}
(i_1, \ldots, i_k) \DI (T_\lambda, P_\lambda).
\end{equation}
Figure \ref{fig:DI-example}  shows the calculations which give
$$
(2,4,3) \DI
\left(
{\beginpicture
\setcoordinatesystem units <0.38cm,0.38cm> % sets scale
\setplotarea x from 9.5 to 12, y from -1 to 1    % sets plot size up
\linethickness=0.5pt      
\putrectangle corners at  9.5 1.5 and 10.5 0.5
\putrectangle corners at 10.5 1.5 and 11.5 0.5
\putrectangle corners at 11.5 1.5 and 12.5 0.5
\putrectangle corners at 12.5 1.5 and 13.5 0.5
\putrectangle corners at  9.5 0.5 and 10.5 -.5
\putrectangle corners at  9.5 -.5 and 10.5 -1.5
\put{1} at 10 1 \put{2} at 11 1  \put{3} at 12 1 \put{6} at 13 1   
\put{4} at 10 0 
\put{5} at 10 -1
\endpicture},
\left(
{\beginpicture
\setcoordinatesystem units <0.2cm,0.2cm> % sets scale
\setplotarea x from 0 to 0, y from -1 to 1    % sets plot size up
\linethickness=0.5pt      
\putrectangle corners at 0 1 and 1 0
\putrectangle corners at 1 1 and 2 0
\putrectangle corners at 2 1 and 3 0
\putrectangle corners at 3 1 and 4 0
\putrectangle corners at 4 1 and 5 0
\putrectangle corners at 5 1 and 6 0
\endpicture},
{\beginpicture
\setcoordinatesystem units <0.2cm,0.2cm> % sets scale
\setplotarea x from 0 to 0, y from -1 to 1    % sets plot size up
\linethickness=0.5pt      
\putrectangle corners at 0 1 and 1 0
\putrectangle corners at 1 1 and 2 0
\putrectangle corners at 2 1 and 3 0
\putrectangle corners at 3 1 and 4 0
\putrectangle corners at 4 1 and 5 0
\endpicture},
{\beginpicture
\setcoordinatesystem units <0.2cm,0.2cm> % sets scale
\setplotarea x from 0 to 0, y from -1 to 1    % sets plot size up
\linethickness=0.5pt      
\putrectangle corners at 0 1 and 1 0
\putrectangle corners at 1 1 and 2 0
\putrectangle corners at 2 1 and 3 0
\putrectangle corners at 3 1 and 4 0
\putrectangle corners at 4 1 and 5 0
\putrectangle corners at 0 0 and 1 -1
\endpicture},
{\beginpicture
\setcoordinatesystem units <0.2cm,0.2cm> % sets scale
\setplotarea x from 0 to 0, y from -1 to 1    % sets plot size up
\linethickness=0.5pt      
\putrectangle corners at 0 1 and 1 0
\putrectangle corners at 1 1 and 2 0
\putrectangle corners at 2 1 and 3 0
\putrectangle corners at 3 1 and 4 0
\putrectangle corners at 0 0 and 1 -1
\endpicture},
{\beginpicture
\setcoordinatesystem units <0.2cm,0.2cm> % sets scale
\setplotarea x from 0 to 0, y from -1 to 1    % sets plot size up
\linethickness=0.5pt      
\putrectangle corners at 0 1 and 1 0
\putrectangle corners at 1 1 and 2 0
\putrectangle corners at 2 1 and 3 0
\putrectangle corners at 3 1 and 4 0
\putrectangle corners at 0 0 and 1 -1
\putrectangle corners at 1 0 and 2 -1
\endpicture},
{\beginpicture
\setcoordinatesystem units <0.2cm,0.2cm> % sets scale
\setplotarea x from 0 to 0, y from -1 to 1    % sets plot size up
\linethickness=0.5pt      
\putrectangle corners at 0 1 and 1 0
\putrectangle corners at 1 1 and 2 0
\putrectangle corners at 2 1 and 3 0
\putrectangle corners at 3 1 and 4 0
\putrectangle corners at 0 0 and 1 -1
\endpicture},
{\beginpicture
\setcoordinatesystem units <0.2cm,0.2cm> % sets scale
\setplotarea x from 0 to 0, y from -1 to 1    % sets plot size up
\linethickness=0.5pt      
\putrectangle corners at 0 1 and 1 0
\putrectangle corners at 1 1 and 2 0
\putrectangle corners at 2 1 and 3 0
\putrectangle corners at 3 1 and 4 0
\putrectangle corners at 0 0 and 1 -1
\putrectangle corners at 0 -1 and 1 -2
\endpicture}
\right)
\right).
$$

\begin{figure}
\label{fig:DI-example}
\caption{Delete-Insertion of $(2,4,3)$}
$$
 {\beginpicture
\setcoordinatesystem units <0.37cm,0.37cm> % sets scale
\setplotarea x from 0 to 13, y from -16 to 3    % sets plot size up
\linethickness=0.5pt                        % sets line thickness

\put{$j$} at 0 6
\put{$i_j$} at 5 6
\put{$T^{(j)}$} at 10 6

\plot -1 5 16 5 /

\put{$0$}    at 0  3    
\put{${1\over2}$}  at 0  0    \put{$2~\jdt$} at 6 0
\put{$1 $}   at 0 -3   \put{$2~\RS$} at 6 -3
\put{$1{1\over2}$}  at 0 -6   \put{$4~\jdt$} at 6 -6
\put{$2 $}   at 0 -9   \put{$4 ~\RS$} at 6 -9
\put{$2{1\over2}$}  at 0 -12   \put{$3~\jdt$} at 6 -12
\put{$3 $}   at 0 -15  \put{$3~\RS$} at 6 -15

\putrectangle corners at 9.5 2.5 and 10.5 3.5
\putrectangle corners at 10.5 2.5 and 11.5 3.5
\putrectangle corners at 11.5 2.5 and 12.5 3.5
\putrectangle corners at 12.5 2.5 and 13.5 3.5
\putrectangle corners at 13.5 2.5 and 14.5 3.5
\putrectangle corners at 14.5 2.5 and 15.5 3.5
\put{1} at 10 3 \put{2} at 11 3 \put{3} at 12 3 \put{4} at 13 3 \put{5} at 14 3 \put{6} at 15 3

\putrectangle corners at  9.5 -0.5 and 10.5 0.5
\putrectangle corners at 10.5 -0.5 and 11.5 0.5
\putrectangle corners at 11.5 -0.5 and 12.5 0.5
\putrectangle corners at 12.5 -0.5 and 13.5 0.5
\putrectangle corners at 13.5 -0.5 and 14.5 0.5
\put{1} at 10 0 \put{3} at 11 0 \put{4} at 12 0 \put{5} at 13 0 \put{6} at 14 0 

\putrectangle corners at  9.5 -2.5 and 10.5 -3.5
\putrectangle corners at 10.5 -2.5 and 11.5 -3.5
\putrectangle corners at 11.5 -2.5 and 12.5 -3.5
\putrectangle corners at 12.5 -2.5 and 13.5 -3.5
\putrectangle corners at 13.5 -2.5 and 14.5 -3.5
\putrectangle corners at  9.5 -3.5 and 10.5 -4.5
\put{1} at 10 -3 \put{2} at 11 -3  \put{4} at 12 -3 \put{5} at 13 -3  \put{6} at 14 -3 
\put{3} at 10 -4

\putrectangle corners at  9.5 -5.5 and 10.5 -6.5
\putrectangle corners at 10.5 -5.5 and 11.5 -6.5
\putrectangle corners at 11.5 -5.5 and 12.5 -6.5
\putrectangle corners at 12.5 -5.5 and 13.5 -6.5
\putrectangle corners at  9.5 -6.5 and 10.5 -7.5
\put{1} at 10 -6 \put{2} at 11 -6  \put{5} at 12 -6 \put{6} at 13 -6   
\put{3} at 10 -7
 
\putrectangle corners at  9.5 -8.5 and 10.5 -9.5
\putrectangle corners at 10.5 -8.5 and 11.5 -9.5
\putrectangle corners at 11.5 -8.5 and 12.5 -9.5
\putrectangle corners at 12.5 -8.5 and 13.5 -9.5
\putrectangle corners at  9.5 -9.5 and 10.5 -10.5
\putrectangle corners at  10.5 -9.5 and 11.5 -10.5
\put{1} at 10 -9 \put{2} at 11 -9  \put{4} at 12 -9 \put{6} at 13 -9  
\put{3} at 10 -10 \put{5} at 11 -10

\putrectangle corners at  9.5 -11.5 and 10.5 -12.5
\putrectangle corners at 10.5 -11.5 and 11.5 -12.5
\putrectangle corners at 11.5 -11.5 and 12.5 -12.5
\putrectangle corners at  12.5 -11.5 and 13.5 -12.5
\putrectangle corners at  9.5 -12.5 and 10.5 -13.5
\put{1} at 10 -12 \put{2} at 11 -12  \put{4} at 12 -12 \put{6} at 13 -12 
\put{5} at 10 -13

\putrectangle corners at  9.5 -14.5 and 10.5 -15.5
\putrectangle corners at 10.5 -14.5 and 11.5 -15.5
\putrectangle corners at 11.5 -14.5 and 12.5 -15.5
\putrectangle corners at  12.5 -14.5 and 13.5 -15.5
\putrectangle corners at  9.5 -15.5 and 10.5 -16.5
\putrectangle corners at 9.5 -16.5 and 10.5 -17.5
\put{1} at 10 -15 \put{2} at 11 -15  \put{3} at 12 -15 \put{6} at 13 -15 
\put{4} at 10 -16 
\put{5} at 10 -17

\endpicture}
$$
\end{figure}

\begin{thm} The function $(i_1, \ldots, i_k) \DI (T_\lambda,P_\lambda)$ provides a bijection
between sequences in $\left\{\, (i_1, \ldots, i_k) \, \Big|\, 1 \le i_j \le n\,\right\}$ and
$\bigsqcup_{\lambda \in \Lambda_n^k} \SYT(\lambda) \times \PT_k(\lambda)$ and thus gives a 
combinatorial proof of (\ref{eq:identity1}).
\end{thm}

\begin{proof}  We prove the theorem by constructing the inverse of $\DI$.  

Let $\lambda^{(j+{1\over2})} \subseteq \lambda^{(j+1)}$ with
$\lambda^{(j+1)} \in
\Lambda_n^{(j+1)}$ and $\lambda^{(j+\half)} \in \Lambda_{n-1}^{j}$, and  let  $T^{(j+1)}$ be a  standard  
tableau of shape $\lambda^{(j+1)}$. We can uniquely determine $i_{j+1}$ and a  tableau $T^{(j+{1\over2})}$ of shape
$\lambda^{(j+{1\over2})}$ such that $T^{(j+1)} = (i_{j+1} \RS T^{(j+{1\over2} )}).$  To do this let $b$ be the box
in $\lambda^{(j+1)}/\lambda^{(j+{1\over2})}$. Then let $i_{j+1}$ and $T^{(j+{1\over2})}$ be the result of uninserting the number 
in box $b$ of $T^{(j+1)}$  (using the fact that RSK insertion is invertible).

Now let $T^{(j+{1\over2})}$ be a tableau of shape $\lambda^{(j+{1\over2})} \in \Lambda_{n-1}^{(j)}$ with  increasing rows and
columns and entries $\{1, \ldots, n\} \setminus \{i_{j+1}\}$ and let $\lambda^{(j)} \subseteq \lambda^{(j+{1\over2})}$
with
$\lambda^{(j)} \in \Lambda^j_n$. We can uniquely produce a standard  tableau $T^{(j)}$ such that
$T^{(j+{1\over2})} = (i_{j+1} \jdt T^{(j)})$.  To do this, let $b$ be the box in $\lambda^{(j)}/\lambda^{(j+{1\over2})}$, put 
$i_{j+1}$ in position $b$ of $T^{(j+{1\over2})}$, and perform the inverse of jeu de taquin to produce $T^{(j)}$.  That is,  move $i_{j+1}$ into a standard position by iteratively swapping it with the larger of the numbers just above
it or just to its left.

Given $\lambda \in \Lambda_n^k$ and $(P_\lambda,T_\lambda) \in \SYT(\lambda) \times \PT_k(\lambda)$, 
we apply the process above to $\lambda^{(k-\half)} \subseteq \lambda^{(k)}$ and
$T^{(k)} = T_\lambda$  producing $i_k$ and $T^{(k-1)}$. 
Continuing this way, we can produce $i_k, i_{k-1}, \ldots, i_1$ and $T^{(k)}, T^{(k-1)}, \ldots, T^{(1)}$ such  
that $(i_1, \ldots, i_k) \DI (T_\lambda,P_\lambda).$
\end{proof}

\section{A Bijective Proof of $B(2k) = \displaystyle{\sum_{\lambda \in
\Gamma_k}} (m^\lambda_k)^2$}

If we view $\CC A_k(n)$ as a bimodule for $\CC A_k(n) \otimes \CC A_k(n)$, with the first
tensor factor acting by left multiplication on $\CC A_k(n)$ and the second tensor factor acting
by right multiplication on $\CC A_k(n)$, then the decomposition into irreducibles for
$\CC A_k(n) \otimes \CC A_k(n)$ is
$
\CC A_k(n) \cong \bigoplus_{\lambda\in\Gamma_k} M^\lambda_k \otimes \bar{M}^\lambda_k,
$
where $\bar{M}^\lambda_k$ is the irreducible right $\CC A_k(n)$-module indexed by $\lambda$
(so $M^\lambda_k \cong \bar{M}^\lambda_k$).  This is a generalization of the decomposition 
of the group algebra of a finite group.
Comparing dimensions on both sides gives,
\begin{equation}
B(2k) = \sum_{\lambda\in\Gamma_k} (m^\lambda_k)^2.
\label{eq:identity2}
\end{equation}
For example using the numbers $m^\lambda_3$ from the bottom row of Figure 2, we have
$
B(6) = 203 = 5^2 + 10^2 + 6^2 + 6^2 + 1^2 + 2^2 + 1^2.
$

To give a combinatorial proof of (\ref{eq:identity2}), we find a bijection
of the form
\begin{equation}
A_{k} \longleftrightarrow 
\bigsqcup_{\lambda \in \Gamma_k} \PT_k(\lambda) \times \PT_k(\lambda),
\end{equation}
by constructing a function that takes a set partition $d \in A_{k}$ and produces
a pair $(P_\lambda,Q_\lambda)$ of vacillating tableaux.
In Section 5 we show that our bijection restricts to work for all the diagram subalgebras described
in Section 2.
In particular, if $d$ corresponds to a permutation in $S_k$ then our bijection is the usual RSK
algorithm described in Section 3.2.

We will draw diagrams $d \in A_k$ using a {\it standard representation} as single row with the vertices in order $1, \ldots, 2k$, where we relabel vertex $j'$ with with the label $2k - j + 1$. We draw the edges of the standard representation of $d \in A_k$ in a specific way:
connect vertices $i$ and $j$, with $i \le j$, if and only if $i$ and $j$ are related in $d$
and there does not exist $k$ related to $i$ and $j$ with $i < k < j$.  In this way, each vertex is connected  only to its nearest neighbors in its block.  
For example,
$$
 
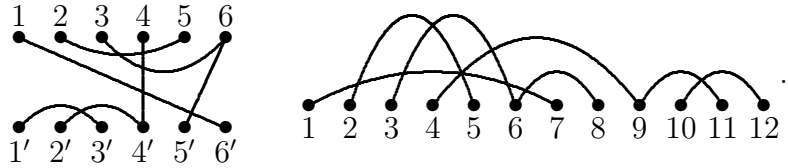
{\beginpicture
\setcoordinatesystem units <0.55cm,0.3cm> % sets scale
\setplotarea x from 0 to 7, y from -3 to 3    % sets plot size up
\linethickness=0.5pt                        % sets line thickness
\put{1} at 1 3
\put{2} at 2 3
\put{3} at 3 3
\put{4} at 4 3
\put{5} at 5 3
\put{6} at 6 3
\put{$1'$} at  1 -3
\put{$2'$} at 2 -3
\put{$3'$} at 3 -3
\put{$4'$} at 4 -3
\put{$5'$} at 5 -3
\put{$6'$} at 6 -3
\put{$\bullet$} at 1 -2 \put{$\bullet$} at 1 2
\put{$\bullet$} at 2 -2 \put{$\bullet$} at 2 2
\put{$\bullet$} at 3 -2 \put{$\bullet$} at 3 2
\put{$\bullet$} at 4 -2 \put{$\bullet$} at 4 2
\put{$\bullet$} at 5 -2 \put{$\bullet$} at 5 2
\put{$\bullet$} at 6 -2 \put{$\bullet$} at 6 2
\plot 1 2 6 -2  /
\plot 4 2 4 -2 /
\plot 6 2 5 -2 /
\setquadratic
\plot 2 2 3.5 1.25 5 2 /
\plot 3 2 4.5 .5 6 2 /
\plot 1 -2 2 -1 3 -2 /
\plot 2 -2 3 -1 4 -2 /
\endpicture}
 {\beginpicture
\setcoordinatesystem units <.55cm,0.3cm> % sets scale
\setplotarea x from 0 to 12, y from -2 to 3    % sets plot size up
\linethickness=0.5pt
\put{1} at 1 -2
\put{2} at 2 -2
\put{3} at 3 -2
\put{4} at 4 -2
\put{5} at 5 -2
\put{6} at 6 -2
\put{7} at 7 -2
\put{8} at 8 -2
\put{9} at 9 -2
\put{10} at 10 -2
\put{11} at 11 -2
\put{12} at 12 -2
\put{$\bullet$} at 1 -1 \put{$\bullet$} at 7 -1
\put{$\bullet$} at 2 -1 \put{$\bullet$} at 8 -1
\put{$\bullet$} at 3 -1 \put{$\bullet$} at 9 -1
\put{$\bullet$} at 4 -1 \put{$\bullet$} at 10 -1
\put{$\bullet$} at 5 -1 \put{$\bullet$} at 11 -1
\put{$\bullet$} at 6 -1 \put{$\bullet$} at 12 -1
\setquadratic
\plot 10 -1 11 0.5 12 -1 /
\plot 9 -1  10 0.5 11 -1 /
\plot 4 -1 6.5 2 9 -1 /
\plot 6 -1   7 0.5 8 -1 /
\plot 1 -1   4 0.5 7 -1 /
\plot 3 -1 4.5 3.0 6 -1 /
\plot 2 -1 3.5 3.0 5 -1 /
\endpicture}.
$$

We label each edge $e$ of the diagram $d$ with  $2k+1- \ell$, where $\ell$ is the right vertex of $e$.  Define the {\it insertion sequence\/} of a  diagram to be the sequence $E = (E_j)$
indexed by $j$ in the sequence $ {1\over2}, 1, 1{1\over2}, \ldots, 2k-1, 2k-{1\over2}, 2k$, where
\begin{eqnarray*}
E_{j} &=&
\begin{cases}
a, & \hbox{if vertex $j$ is the left endpoint of edge $a$}, \\
\emptyset, & \hbox{if vertex $j$ is not a left endpoint}, \\
\end{cases}  \\
E_{j-{1\over2}} &=&
\begin{cases}
a, & \hbox{if vertex $j$ is the right endpoint of edge $a$}, \\
\emptyset, & \hbox{if vertex $j$ is not a right endpoint}. \\
\end{cases} \\
\end{eqnarray*}

\begin{xmp}  The edge labeling for the diagram of the set partition 
$\big\{ \{1,3,4'\}, \{2,1'\},$ $\{4,3',2'\}\}$ is given by
$$
 
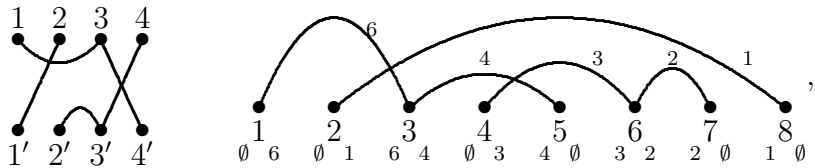
{\beginpicture
\setcoordinatesystem units <0.55cm,0.3cm> % sets scale
\setplotarea x from 0 to 5, y from -3 to 3    % sets plot size up
\linethickness=0.5pt                        % sets line thickness
\put{1} at 1 3
\put{2} at 2 3
\put{3} at 3 3
\put{4} at 4 3
\put{$1'$} at  1 -3
\put{$2'$} at 2 -3
\put{$3'$} at 3 -3
\put{$4'$} at 4 -3
\put{$\bullet$} at 1 -2 \put{$\bullet$} at 1 2
\put{$\bullet$} at 2 -2 \put{$\bullet$} at 2 2
\put{$\bullet$} at 3 -2 \put{$\bullet$} at 3 2
\put{$\bullet$} at 4 -2 \put{$\bullet$} at 4 2
\plot 2 2 1 -2  /
\plot 3 2 4 -2 /
\plot 4 2 3 -2 /
\setquadratic
\plot 1 2 2 1 3 2 /
\plot 2 -2 2.5 -1 3 -2 /
\endpicture}
 {\beginpicture
\setcoordinatesystem units <1cm,0.3cm> % sets scale
\setplotarea x from 0 to 8, y from -2 to 3    % sets plot size up
\linethickness=0.5pt
\put{1} at 1 -2
\put{2} at 2 -2
\put{3} at 3 -2
\put{4} at 4 -2
\put{$5$} at 5 -2
\put{$6$} at 6 -2
\put{$7$} at 7 -2
\put{$8$} at 8 -2
\put{$\bullet$} at 1 -1 \put{$\bullet$} at 7 -1
\put{$\bullet$} at 2 -1 \put{$\bullet$} at 8 -1
\put{$\bullet$} at 3 -1 
\put{$\bullet$} at 4 -1 
\put{$\bullet$} at 5 -1 
\put{$\bullet$} at 6 -1
\setquadratic
\plot 1 -1   2.0 3.0 3 -1 /
\plot 4 -1   5 1 6 -1 /
\plot 6 -1   6.5 .75 7 -1 /
\plot 2 -1   5 3 8 -1 /
\plot 3 -1   4 .5 5 -1 /
\put{$\scriptstyle{1}$} at 7.5 1.2
\put{$\scriptstyle{2}$} at 6.5 1.2
\put{$\scriptstyle{3}$} at 5.5 1.2
\put{$\scriptstyle{4}$} at 4 1.2
\put{$\scriptstyle{6}$} at 2.5 2.5
\put{$\scriptstyle{\emptyset}$} at .8 -3
\put{$\scriptstyle{6}$} at 1.2 -3
\put{$\scriptstyle{\emptyset}$} at 1.8 -3
\put{$\scriptstyle{1}$} at 2.2 -3
\put{$\scriptstyle{6}$} at 2.8 -3
\put{$\scriptstyle{4}$} at 3.2 -3
\put{$\scriptstyle{\emptyset}$} at 3.8 -3
\put{$\scriptstyle{3}$} at 4.2 -3
\put{$\scriptstyle{4}$} at 4.8 -3
\put{$\scriptstyle{\emptyset}$} at 5.2 -3
\put{$\scriptstyle{3}$} at 5.8 -3
\put{$\scriptstyle{2}$} at 6.2 -3
\put{$\scriptstyle{2}$} at 6.8 -3
\put{$\scriptstyle{\emptyset}$} at 7.2 -3
\put{$\scriptstyle{1}$} at 7.8 -3
\put{$\scriptstyle{\emptyset}$} at 8.2 -3
\endpicture},
$$
and the insertion sequence is
$$
\begin{tabular}{c|ccccccccccccccccc}
\hline
$j$ & $\half$ & $1$ & $1\half$ & $2$ & $2\half$ & $3$ & $3\half$
& $4$  & $4\half$ & $5$ & $5\half$ & $6$ & $6\half$ & $7$ & $7\half$ & $8$ \\ \hline
$E_j$ & $\emptyset$ & 6 & $\emptyset$ & 1 & 6 & 4 &  $\emptyset$  & 3 & 4 & $\emptyset$  & 3 & 2 & 2 & $\emptyset$  & 1 & $\emptyset$ \\
\hline
\end{tabular}.
$$
\label{example:insertionsequence}
\end{xmp}

The insertion sequence of a standard diagram completely determines the edges, and thus the connected components, of the diagram, so the following proposition follows immediately.

\begin{prop}  $d \in A_k$ is completely determined by its insertion sequence.
\end{prop}

For $d \in A_k$ with insertion sequence $E = (E_j)$, we will produce a pair $(P_\lambda, Q_\lambda)$ of vacillating tableaux. Begin with the empty tableaux,
\begin{equation}
T^{(0)} =  \emptyset.
\end{equation}
Then recursively define standard  tableaux  $T^{(j+{1\over2})}$ and $T^{(j+1)}$ by
\begin{equation}
\begin{array}{rcl}
T^{(j+{1\over2})} &=& 
\begin{cases}
 E_{j+{1\over2}} \jdt T^{(j)}, & \hbox{if $E_{j+{1\over2}} \not= \emptyset$}, \\
T^{(j)}, & \hbox{if $E_{j+{1\over2}} = \emptyset$,} \\
\end{cases} \\
T^{(j+1)} &=& 
\begin{cases}
E_{j+1} \RS T^{(j+{1\over2})}, & \hbox{if $E_{j+1} \not= \emptyset$}, \\
T^{(j+{1\over2})}, & \hbox{if $E_{j+1} = \emptyset$,} \\
\end{cases} \\
\end{array} \qquad\quad 0 \le j \le 2k-1.
\end{equation}
Let $\lambda^{(i)}$ be the shape of $T^{(i)}$, let $\lambda^{(i+{1\over2})}$ be the shape
of $T^{(i+{1\over2})}$,  and let $\lambda = \lambda^{(k)}$. Define
\begin{eqnarray*}
Q_\lambda &=& \left(\emptyset, \lambda^{({1\over2})}, \lambda^{(1)}, \ldots,\lambda^{(k-{1\over2})}, \lambda^{(k)}
\right) \in \PT(\lambda),
\\ P_\lambda &=& \left( \lambda^{(2k)}, \lambda^{(2k-{1\over2})}, \ldots, \lambda^{(k+{1\over2})},\lambda^{(k)}
\right) \in \PT(\lambda).
\end{eqnarray*}
In this way, we associate a pair of vacillating tableaux $(P_\lambda, Q_\lambda)$ to a set partition
$d \in A_{k}$, which we denote by
$$
d  \to (P_\lambda,Q_\lambda).
$$
In Figure \ref{fig:insertdiag}, we see the insertion of the diagram in Example \ref{example:insertionsequence}.  The result of this is to pair the diagram in Example \ref{example:insertionsequence}
 with the following pair of vacillating tableaux,
\begin{eqnarray*}
Q_\lambda &=& \left( 
\emptyset, ~\emptyset,
{\beginpicture
\setcoordinatesystem units <0.3cm,0.3cm>         % sets scale
\setplotarea x from .8 to 2.2, y from -1 to 1   % sets plot size up
\linethickness=0.5pt   
\putrectangle corners at 1  0 and  2 1
\endpicture},
{\beginpicture
\setcoordinatesystem units <0.3cm,0.3cm>         % sets scale
\setplotarea x from .8 to 2.2, y from -1 to 1   % sets plot size up
\linethickness=0.5pt   
\putrectangle corners at 1  0 and  2 1
\endpicture},
{\beginpicture
\setcoordinatesystem units <0.3cm,0.3cm>         % sets scale
\setplotarea x from .8 to 2.2, y from -1 to 1   % sets plot size up
\linethickness=0.5pt   
\putrectangle corners at 1  0 and  2 1
\putrectangle corners at 1 -1 and 2 0
\endpicture},
{\beginpicture
\setcoordinatesystem units <0.3cm,0.3cm>         % sets scale
\setplotarea x from .8 to 2.2, y from -1 to 1   % sets plot size up
\linethickness=0.5pt   
\putrectangle corners at 1  0 and  2 1
\endpicture},
{\beginpicture
\setcoordinatesystem units <0.3cm,0.3cm>         % sets scale
\setplotarea x from .8 to 3.2, y from -1 to 1   % sets plot size up
\linethickness=0.5pt   
\putrectangle corners at 1  0 and  2 1
\putrectangle corners at 2  0 and  3 1
\endpicture},
{\beginpicture
\setcoordinatesystem units <0.3cm,0.3cm>         % sets scale
\setplotarea x from .8 to 3.2, y from -1 to 1   % sets plot size up
\linethickness=0.5pt   
\putrectangle corners at 1  0 and  2 1
\putrectangle corners at 2  0 and  3 1
\endpicture},
{\beginpicture
\setcoordinatesystem units <0.3cm,0.3cm>         % sets scale
\setplotarea x from .8 to 3.2, y from -1 to 1   % sets plot size up
\linethickness=0.5pt   
\putrectangle corners at 1  0 and  2 1
\putrectangle corners at 2  0 and  3 1
\putrectangle corners at 1 -1 and 2 0
\endpicture}
\right), \\
P_\lambda &=& \left( 
\emptyset, ~\emptyset,
{\beginpicture
\setcoordinatesystem units <0.3cm,0.3cm>         % sets scale
\setplotarea x from .8 to 2.2, y from -1 to 1   % sets plot size up
\linethickness=0.5pt   
\putrectangle corners at 1  0 and  2 1
\endpicture},
{\beginpicture
\setcoordinatesystem units <0.3cm,0.3cm>         % sets scale
\setplotarea x from .8 to 2.2, y from -1 to 1   % sets plot size up
\linethickness=0.5pt   
\putrectangle corners at 1  0 and  2 1
\endpicture},
{\beginpicture
\setcoordinatesystem units <0.3cm,0.3cm>         % sets scale
\setplotarea x from .8 to 3.2, y from -1 to 1   % sets plot size up
\linethickness=0.5pt   
\putrectangle corners at 1  0 and  2 1
\putrectangle corners at 2 0 and 3 1
\endpicture},
{\beginpicture
\setcoordinatesystem units <0.3cm,0.3cm>         % sets scale
\setplotarea x from .8 to 2.2, y from -1 to 1   % sets plot size up
\linethickness=0.5pt   
\putrectangle corners at 1  0 and  2 1
\endpicture},
{\beginpicture
\setcoordinatesystem units <0.3cm,0.3cm>         % sets scale
\setplotarea x from .8 to 2.2, y from -1 to 1   % sets plot size up
\linethickness=0.5pt   
\putrectangle corners at 1 0  and  2 1
\putrectangle corners at 2 0 and 3 1
\endpicture},
{\beginpicture
\setcoordinatesystem units <0.3cm,0.3cm>         % sets scale
\setplotarea x from .8 to 2.2, y from -1 to 1   % sets plot size up
\linethickness=0.5pt   
\putrectangle corners at 1 0  and  2 1
\putrectangle corners at 2 0 and 3 1
\endpicture},
{\beginpicture
\setcoordinatesystem units <0.3cm,0.3cm>         % sets scale
\setplotarea x from .8 to 3.2, y from -1 to 1   % sets plot size up
\linethickness=0.5pt   
\putrectangle corners at 1  0 and  2 1
\putrectangle corners at 2  0 and  3 1
\putrectangle corners at 1 -1 and 2 0
\endpicture}
\right).
\end{eqnarray*}

\begin{figure} \label{fig:insertdiag}
\caption{Insertion of the Set Partition in Example \ref{example:insertionsequence}. }

$$
\begin{tabular}{c|ccccccccccccccccc}
\hline
$j$ & $\half$ & $1$ & $1\half$ & $2$ & $2\half$ & $3$ & $3\half$
& $4$  & $4\half$ & $5$ & $5\half$ & $6$ & $6\half$ & $7$ & $7\half$ & $8$ \\ \hline
$E_j$ & $\emptyset$ & 6 & $\emptyset$ & 1 & 6 & 4 &  $\emptyset$  & 3 & 4 & $\emptyset$  & 3 & 2 & 2 & $\emptyset$  & 1 & $\emptyset$ \\
\hline
\end{tabular}
$$

$$
\begin{tabular}{cccc}
$j$ & $E_j$ & & $T^{(j)}$ \\ \hline \\
$0$ & & & $\emptyset$ \\ \\
${1 \over 2}$ & $\emptyset$&$ \jdt$ & $\emptyset$ \\ \\
$1$ & $6$&$\RS $ & 
{\beginpicture
\setcoordinatesystem units <0.4cm,0.4cm>         
\putrectangle corners at 1  -.2 and  2 .8 \put{6} at 1.5 .3 \endpicture} \\\\
$1{1 \over 2}$ & $\emptyset$&$\jdt $ & 
{\beginpicture
\setcoordinatesystem units  <0.4cm,0.4cm>       
\putrectangle corners at 1  -.2 and  2 .8 \put{6} at 1.5 .3 \endpicture} \\\\
$2$ & $1$&$\RS $ & 
{\beginpicture
\setcoordinatesystem units <0.4cm,0.4cm>          
\putrectangle corners at 1  -.2 and  2 .8 \put{1} at 1.5 .3 
\putrectangle corners at 1  -1.2 and  2 -.2 \put{6} at 1.5 -.7 \endpicture} \\\\
$2{1 \over 2}$ & $6$&$\jdt $ & 
{\beginpicture
\setcoordinatesystem units <0.4cm,0.4cm>         
\putrectangle corners at 1  -.2 and  2 .8 \put{1} at 1.5 .3 \endpicture} \\\\
$3$ & $4$&$\RS $ & 
{\beginpicture
\setcoordinatesystem units <0.4cm,0.4cm>          
\putrectangle corners at 1  -.2 and  2 .8 \put{1} at 1.5 .3 
\putrectangle corners at 2  -.2 and  3 .8 \put{4} at 2.5 .3  \endpicture} \\\\
$3{1 \over 2}$ & $\emptyset$&$\jdt $ & 
{\beginpicture
\setcoordinatesystem units <0.4cm,0.4cm>          
\putrectangle corners at 1  -.2 and  2 .8 \put{1} at 1.5 .3 
\putrectangle corners at 2  -.2 and  3 .8 \put{4} at 2.5 .3  \endpicture}  \\\\
$4$ & $3$&$\RS $ & 
{\beginpicture
\setcoordinatesystem units <0.4cm,0.4cm>          
\putrectangle corners at 1  -.2 and  2 .8 \put{1} at 1.5 .3 
\putrectangle corners at 1  -1.2 and  2 -.2 \put{4} at 1.5 -.7
\putrectangle corners at 2  -.2 and  3 .8 \put{3} at 2.5 .3  \endpicture} \\\\
\end{tabular}
\hskip1truein
\begin{tabular}{cccc}
$j$ & $E_j$ & & $T^{(j)}$ \\ \hline \\
$8$ & $\emptyset$&$ \RS$ & $\emptyset$ \\ \\
$7{1 \over 2}$ & $1$&$ \jdt$ & $\emptyset$ \\ \\
$7 $ & $\emptyset$&$\RS $ & 
{\beginpicture
\setcoordinatesystem units <0.4cm,0.4cm>         
\putrectangle corners at 1  -.2 and  2 .8 \put{1} at 1.5 .3 \endpicture} \\\\
$6{1 \over 2}$ & $2$&$\jdt $ & 
{\beginpicture
\setcoordinatesystem units <0.4cm,0.4cm> 
\putrectangle corners at 1  -.2 and  2 .8 \put{1} at 1.5 .3 \endpicture} \\\\
$6$ & $2$&$\RS $ & 
{\beginpicture
\setcoordinatesystem units <0.4cm,0.4cm> 
\putrectangle corners at 1  -.2 and  2 .8 \put{1} at 1.5 .3 
\putrectangle corners at 2  -.2 and  3 .8 \put{2} at 2.5 .3  \endpicture} \\\\
$5{1 \over 2}$ & $3$&$\jdt $ & 
{\beginpicture
\setcoordinatesystem units <0.4cm,0.4cm> 
\putrectangle corners at 1  -.2 and  2 .8 \put{1} at 1.5 .3 \endpicture} \\\\
$5$ & $\emptyset$&$\RS $ & 
{\beginpicture
\setcoordinatesystem units <0.4cm,0.4cm> 
\putrectangle corners at 1  -.2 and  2 .8 \put{1} at 1.5 .3 
\putrectangle corners at 2  -.2 and  3 .8 \put{3} at 2.5 .3  \endpicture} \\\\
$4{1 \over 2}$ & 4&$\jdt $ & 
{\beginpicture
\setcoordinatesystem units <0.4cm,0.4cm> 
\putrectangle corners at 1  -.2 and  2 .8 \put{1} at 1.5 .3 
\putrectangle corners at 2  -.2 and  3 .8 \put{3} at 2.5 .3  \endpicture}  \\\\
$4$ & & & 
{\beginpicture
\setcoordinatesystem units <0.4cm,0.4cm>      
\putrectangle corners at 1  -.2 and  2 .8 \put{1} at 1.5 .3 
\putrectangle corners at 1  -1.2 and  2 -.2 \put{4} at 1.5 -.7
\putrectangle corners at 2  -.2 and  3 .8 \put{3} at 2.5 .3  \endpicture} \\\\
\end{tabular}
$$
\end{figure}

\begin{rem} {\rm We have numbered the edges of the standard diagram of $d$ in increasing order
from right to left, so if $E_{j+{1\over2}} \not= \emptyset$, then $E_{j+{1\over2}}$ will be the largest element of $T^{(j)}$. Thus, in $T^{(j+{1\over2})} = (E_{j+{1\over2}} \jdt T^{(j)})$ we know that $E_{j+{1\over2}}$ is in a corner box, and jeu de taquin simply deletes that box.}
\end{rem}

\begin{thm} The function $\,d \to (P_\lambda,Q_\lambda)$ provides a bijection
between the set partitions $A_k$ and pairs of vacillating tableaux
$\bigsqcup_{\lambda \in \Gamma_k} \PT(\lambda) \times \PT_k(\lambda)$ and thus gives a 
combinatorial proof of (\ref{eq:identity2}).
\end{thm}

\begin{proof}
We prove the theorem by constructing the inverse of $d \to (P_\lambda, Q_\lambda)$.  First, 
we use $Q_\lambda$ followed
by $P_\lambda$ in reverse order to construct the sequence
$\lambda^{({1\over2})},\lambda^{(1)}, \ldots, \lambda^{(2k-{1\over2})},\lambda^{(2k)}$.
We initialize $T^{(2k)} = \emptyset$.

We now show how to construct $T^{(i+{1\over2})}$ and $E_{i+{1}}$ so that  $T^{(i+1)} = (E_{i+{1}} \RS T^{(i+{1\over2})})$.  If 
$\lambda^{(i+{1\over2})} = \lambda^{(i+1)}$, then let $T^{(i+{1\over2})} = T^{(i+1)}$ and $E_{i+{1}} = \emptyset$. Otherwise, $\lambda^{(i+1)}/\lambda^{(i+{1\over2})}$ is a box $b$, and we use uninsertion on the value in $b$ to produce $T^{(i+{1\over2})}$ and $E_{i+{1}}$ such that $T^{(i+1)} = (E_{i+{1}} \RS T^{(i+{1\over2})})$.  Since we uninserted the value in position $b$, we know that $T^{(i+{1\over2})}$ has shape $\lambda^{(i+{1\over2})}$.

We then show how to construct $T^{(i)}$ and $E_{i+{1\over2}}$ so that $T^{(i+{1\over2})} = (E_{i+{1\over2}} \jdt T^{(i)})$.  If $\lambda^{(i)} \lambda^{(i+{1\over2})}$, then let $T^{(i)} = T^{(i+{1\over2})}$ and $E_{i+{1\over2}} = \emptyset$. Otherwise, $\lambda^{(i)}/
\lambda^{(i+{1\over2})}$ is a box $b$.  Let $T^{(i)}$ be the tableau of shape $\lambda^{(i)}$
with the same entries as $T^{(i+1)}$ and having the entry $2k-i$ in box $b$. Let $E_{i+{1\over2}}= 2k -i$.
At any given step $i$ (as we work our way down from step $2k$ to step $1 \over 2$), $2k-i$ is the largest value added to the tableau thus far, so $T^{(i)}$ is standard.  
Furthermore, $T^{(i+{1\over2})} = (E_{i+{1\over2}} \jdt T^{(i)})$, since $E_{i+{1\over2}}=2k-i$ is already in a corner and thus jeu de taquin simply deletes it.

This iterative process will produce $E_{2k}, E_{2k-1}, \ldots, E_{1 \over 2}$ which completely determines
the set partition $d$.  By the way we have constructed $d$, we have $d \to (P_\lambda,Q_\lambda)$.
\end{proof}

\section{Properties of the Insertion Algorithm}

In this section we investigate a number of the nice properties of the insertion algorithm from  Section 4. First, we note that our algorithm restricts to the insertion algorithms for the subalgebras discussed in Section 1: 

\begin{rem} {\rm 
\begin{itemize}
\item[(a)] If $d \in S_k$ is a permutation diagram, then the standard representation of $d$ pairs vertices $1, \ldots k$ with vertices $k+1, \ldots, 2k$. Our algorithm has us insert
on steps $j$ with $1 \le j \le k$, delete on steps $j+{1\over2}$ with $k \le j \le 2k-1$, and do nothing
on all other steps.  If we ignore these ``do-nothing" steps, then our algorithm becomes the usual RSK 
insertion and $P_\lambda$ and $Q_\lambda$ are the insertion and recording tableaux, respectively.

\item[(b)]   If $d \in B_k$ is a Brauer diagram, then each vertex in the standard representation of $d$ is incident
to exactly one edge, so exactly one of $E_{(j-{1\over2})}$ or $E_{j}$ will be nonempty for each $1 \le j \le 2k$. If we ignore the steps where
we  do nothing, then for each vertex we either insert or delete the edge number associated to that vertex. This is precisely the insertion scheme of Sundaram \cite{Sun} (see also \cite{Ter}). This algorithm produces ``oscillating tableaux," which either increase or decrease by one box on each  step.  If we insert all the Brauer diagrams, we produce the Bratteli diagram for $\CC B_k(n)$ (see for example \cite{HR1}).

\item[(c)] If $d \in R_k$ is a rook monoid diagram, then our algorithm has us insert or do nothing on steps $0 <  j \le k$ and delete or do nothing on steps $j + {1 \over 2}$ with $k \le j \le 2k-1$.  If we insert all rook monoid diagrams, we obtain the Bratteli diagram for $\CC R_k$ (see \cite{Ha}) .

\item[(d)] If $d \in A_{k-\frac{1}{2}}$ then the standard representation of $d$ has an edge, labeled $k$, connecting vertices $k$ and $k+1$.  At step $k$ in the insertion process, we will insert this edge into $T^{(k-\half)}$. The label $k$ is larger than all the entries in  $T^{(k-\half)}$, since the entries that are in $T^{(k-\half)}$ correspond to edges whose right endpoint is larger than $k+1$.  Thus we insert this entry at the end of the first row and in the next step, $k + \frac{1}{2}$, we immediately delete it, so that $\lambda^{(k-\half)} = \lambda^{(k+\half)}$.  If we remove $\lambda^{(k)}$ from our pair of paths, we obtain a pair of vacillating tableaux $(P_{\lambda^{(k-\half)}},Q_{\lambda^{(k-\half)}})$ of length $2k-1$, and we get a combinatorial formula for the following dimension identity for $\CC A_{k-\half}(n)$ involving the odd Bell numbers,
\begin{equation}
B(2k-1) = \sum_{{\lambda\in\Gamma_{k-\half}}} (m^\lambda_{k-\half})^2.
\end{equation}
\end{itemize} }
\end{rem}

\begin{prop} If $d \in A_k$ such that $d \to (P_\lambda, Q_\lambda)$, then $|\lambda| = \pn(d)$, where
$\pn(d)$ is the propagating number defined in (\ref{def:propagating}).
\end{prop}

\begin{proof} If we draw a vertical line between vertex $k$ and $k+1$ in the standard representation of $d$, then $\pn(d)$ is the number of edges that cross this line.  The labels of these edges are exactly the numbers which have been inserted but not deleted at step $k$ of our algorithm. Thus the tableau $T^{(k)}$ has exactly these numbers in it and so $\lambda^{(k)}$ consists of $\pn(d)$ boxes.
\end{proof}

Recall that $I_t$ is spanned by the diagrams whose propagating number satisfies $\pn(d) \le t$.   The chain of ideals $\CC I_0(n) \subseteq \CC I_1(n)  \subseteq \cdots \subseteq \CC I_k(n) = \CC A_k(n)$ is a key feature of the ``basic construction" of $\CC A_k(n)$ (see \cite{HR2,Mar3}). The irreducible components of $\CC I_t(n)/(\CC I_{t-1}(n))$ are those indexed by $\lambda \in \Gamma_k$ with $|\lambda| = t$.  The previous proposition tells us that our insertion algorithm respects this basic construction and gives a combinatorial proof of the following dimension identity for $\CC I_t(n)/(\CC I_{t-1}(n))$,
\begin{equation}
\Card\left( \{ d \in A_k \ | \ \pn(d) = t \} \right) = \sum_{{\lambda\in\Gamma_k} \atop \pn(d) = t} (m^\lambda_k)^2.
\end{equation}

\begin{prop} Let  $d \to (P,Q) = (\lambda^{(0)}, \ldots, \lambda^{(2k)}).$  Then $d$ is planar if and only if the length of every partition in the sequence satisfies $\ell(\lambda^{(i)}) \le 1$.
\end{prop}

\begin{proof}  This is a special case of the crossing and nesting theorem in \cite{CDDSY}.  If $d$ is planar, then the 
standard representation of $d$ is also planar. For example,
$$
 
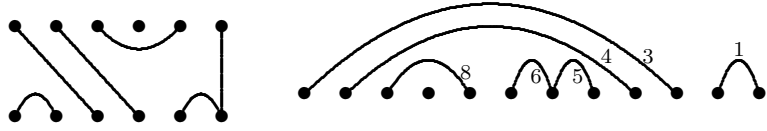
{\beginpicture
\setcoordinatesystem units <0.55cm,0.3cm> % sets scale
\setplotarea x from 0 to 7, y from -2 to 2    % sets plot size up
\linethickness=0.5pt                        % sets line thickness
\put{$\bullet$} at 1 -2 \put{$\bullet$} at 1 2
\put{$\bullet$} at 2 -2 \put{$\bullet$} at 2 2
\put{$\bullet$} at 3 -2 \put{$\bullet$} at 3 2
\put{$\bullet$} at 4 -2 \put{$\bullet$} at 4 2
\put{$\bullet$} at 5 -2 \put{$\bullet$} at 5 2
\put{$\bullet$} at 6 -2 \put{$\bullet$} at 6 2
\plot 1 2 3 -2  /
\plot 2 2 4 -2 /
\plot 6 2 6 -2 /
\setquadratic
\plot 3 2 4 1 5 2 /
\plot 1 -2 1.5 -1 2 -2 /
\plot 5 -2 5.5 -1 6 -2 /
\endpicture}
 {\beginpicture
\setcoordinatesystem units <.55cm,0.3cm> % sets scale
\setplotarea x from 0 to 12, y from -1 to 3    % sets plot size up
\linethickness=0.5pt
\put{$\bullet$} at 1 -1 \put{$\bullet$} at 7 -1
\put{$\bullet$} at 2 -1 \put{$\bullet$} at 8 -1
\put{$\bullet$} at 3 -1 \put{$\bullet$} at 9 -1
\put{$\bullet$} at 4 -1 \put{$\bullet$} at 10 -1
\put{$\bullet$} at 5 -1 \put{$\bullet$} at 11 -1
\put{$\bullet$} at 6 -1 \put{$\bullet$} at 12 -1
\setquadratic
\plot 1 -1 5.5 3 10 -1 /
\plot 2 -1 5.5 2  9 -1 /
\plot 3 -1 4 .5  5 -1 /
\plot 6 -1 6.5 .5  7 -1 /
\plot 7 -1 7.5 .5  8 -1 /
\plot 11 -1 11.5 .5  12 -1 /
\put{$\scriptstyle{1}$} at 11.5 1
\put{$\scriptstyle{3}$} at 9.3 .7
\put{$\scriptstyle{4}$} at 8.3 .7
\put{$\scriptstyle{5}$} at 7.6 -.2
\put{$\scriptstyle{6}$} at 6.6 -.2
\put{$\scriptstyle{8}$} at 4.9 0
\endpicture}.
$$
Let $a < b$ be the labels of two edges in a planar diagram $d$ and suppose that $b$ is inserted before $a$ in the insertion sequence of $d$.   Since $a < b$ the right vertex of $b$ is left of the right vertex of $a$. Since $b$ is inserted before $a$, the left vertex of $b$ is left of the left vertex of $a$.  Since $d$ is planar,  the right vertex of $b$ must then be left of the left vertex of $a$. Thus,  $b$ is inserted and deleted before $a$ is inserted. This means that when we are inserting an edge label $a$ into $T^{(i)}$, we know that $a$ is larger than all the entries of $T^{(i)}$, and so $a$ placed at the end of the first row. Thus at each step of the insertion, we only ever obtain the empty tableau or a tableau consisting of 1 row. Conversely, to obtain a partition that has more than one column, there must be a pair $a < b$ such that  $a$ is inserted after  $b$ is inserted but before $b$ is deleted.  This is necessary for $a$ to bump $b$ into the next row.  By the reasoning above, this forces an edge crossing in the partition diagram.
\end{proof}

\begin{figure}
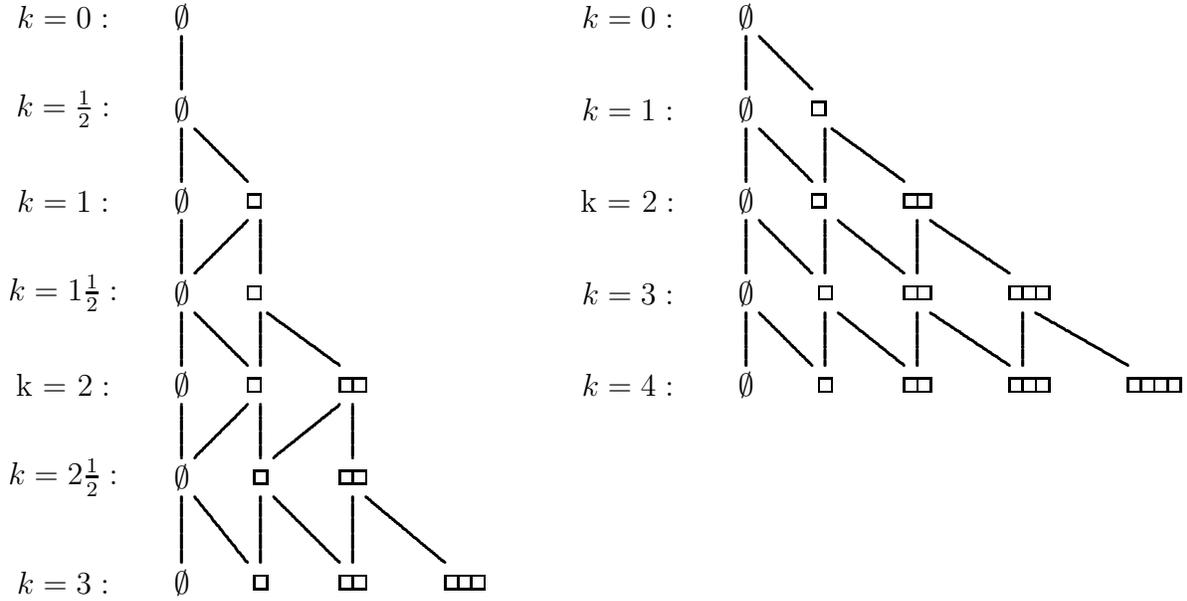

\label{fig:BratteliPlanarPartition}
\caption{Bratteli Diagrams for the Planar Partition Algebra $\CC P_k(n)$ and the Planar Rook Monoid $P\!R_k$}
$$
{\beginpicture
\setcoordinatesystem units <0.175cm,0.175cm>         % sets scale
\setplotarea x from -8 to 26, y from -1 to 26   % sets plot size up
\linethickness=0.5pt                          % sets line thickness
%%%%%%%%%% level 3 %%%%%%%%%%%%%%%%%%%%%%%%%%%%%%%%%%%%%%%%%%%%
\put{$k = 3:$} at  -6 -.5
\put{$\emptyset$} at 3 -.5
\putrectangle corners at 8.5 -1 and 9.5 0
\putrectangle corners at 15 -1 and 16 0
\putrectangle corners at 16 -1 and 17 0
\putrectangle corners at 23 -1 and 24 0
\putrectangle corners at 24 -1 and 25 0
\putrectangle corners at 25 -1 and 26 0
\plot 3 1 3 6  /
\plot 8 1 4 6  /
\plot 9 1 9 6  /
\plot 15 1 10 6  /
\plot 16 1 16 6  /
\plot 23 1 17 6  /
%\plot 21 1 11 6  /
%
%%%%%%%%%% level 2+ %%%%%%%%%%%%%%%%%%%%%%%%%%%%%%%%%%%%%%%%%%%%
\put{$k = 2{1\over2}:$} at  -6 7.5
\put{$\emptyset$} at 3 7.5
\putrectangle corners at 8.5 7  and 9.5 8
\putrectangle corners at 15 7 and 16 8
\putrectangle corners at 16 7 and 17 8
\plot 3 9 3 13  /
\plot 8 13 4 9  /
\plot 9 13 9 9  /
\plot 15 13 10 9  /
\plot 16 13 16 9  /
%%%%%%%%%% level 2 %%%%%%%%%%%%%%%%%%%%%%%%%%%%%%%%%%%%%%%%%%
\put{k = $2:$} at  -6 14.5
\put{$\emptyset$} at 3 14.5
\putrectangle corners at 8  14 and 9   15
\putrectangle corners at 15  14 and 16  15
\putrectangle corners at 16  14 and 17  15
\plot 3 16 3 20  /
\plot 8 16 4 20  /
\plot 9 16 9 20  /
\plot 15 16 9.5 20  /
%%%%%%%%%% level 1+ %%%%%%%%%%%%%%%%%%%%%%%%%%%%%%%%%%%%%%%%%%%%
\put{$k = 1{1\over2}:$} at  -6 21.5
\put{$\emptyset$} at 3 21.5
\putrectangle corners at 8  21 and 9   22
\plot 3 23 3 27  /
\plot 8 27 4 23  /
\plot 9 27 9 23  /
%%%%%%%%%% level 1 %%%%%%%%%%%%%%%%%%%%%%%%%%%%%%%%%%%%%%%%%%%%
\put{$k = 1:$} at  -6 28.5
\put{$\emptyset$} at  3 28.5
\putrectangle corners at 8  28 and 9  29
%
%%%%%%%%%% level 0+ %%%%%%%%%%%%%%%%%%%%%%%%%%%%%%%%%%%%%%%%%%%%
\put{$k = {1\over2}:$} at  -6 35.5
\put{$\emptyset$} at  3 35.5
\plot 3 30 3 34  /
\plot 8 30 4 34  /
%%%%%%%%%% level 0 %%%%%%%%%%%%%%%%%%%%%%%%%%%%%%%%%%%%%%%%%%%%
\put{$k = 0:$} at  -6 42.5
\put{$\emptyset$} at  3 42.5
\plot 3 37 3 41  /
\endpicture} 
\hskip.5truein
{\beginpicture
\setcoordinatesystem units <0.175cm,0.175cm>         % sets scale
\setplotarea x from -8 to 26, y from -1 to 26   % sets plot size up
\linethickness=0.5pt                          % sets line thickness
%%%%%%%%%% level 3 %%%%%%%%%%%%%%%%%%%%%%%%%%%%%%%%%%%%%%%%%%%%
\put{$k = 4:$} at  -6 14.5
\put{$\emptyset$} at 3 14.5
\putrectangle corners at 8.5 14 and 9.5 15
\putrectangle corners at 15 14 and 16 15
\putrectangle corners at 16 14 and 17 15
\putrectangle corners at 23 14 and 24 15
\putrectangle corners at 24 14 and 25 15
\putrectangle corners at 25 14 and 26 15
\plot 3 16 3 20  /
\plot 8 16 4 20  /
\plot 9 16 9 20  /
\plot 15 16 10 20  /
\plot 16 16 16 20  /
\plot 23 16 17 20  /
\plot 24 16 24 20  /
\plot 32 16 25 20  /

\putrectangle corners at 32 14 and 33 15
\putrectangle corners at 33 14 and 34 15
\putrectangle corners at 34 14 and 35 15
\putrectangle corners at 35 14 and 36 15

%%%%%%%%%% level 3 %%%%%%%%%%%%%%%%%%%%%%%%%%%%%%%%%%%%%%%%%%%%
\put{$k = 3:$} at  -6 21.5
\put{$\emptyset$} at 3 21.5
\putrectangle corners at 8.5 21 and 9.5 22
\putrectangle corners at 15 21 and 16 22
\putrectangle corners at 16 21 and 17 22
\putrectangle corners at 23 21 and 24 22
\putrectangle corners at 24 21 and 25 22
\putrectangle corners at 25 21 and 26 22
\plot 3 23 3 27  /
\plot 8 23 4 27  /
\plot 9 23 9 27  /
\plot 15 23 10 27  /
\plot 16 23 16 27  /
\plot 23 23 17 27  /

%%%%%%%%%% level 2 %%%%%%%%%%%%%%%%%%%%%%%%%%%%%%%%%%%%%%%%%%
\put{k = $2:$} at  -6 28.5
\put{$\emptyset$} at 3 28.5
\putrectangle corners at 8  28 and 9   29
\putrectangle corners at 15  28 and 16  29
\putrectangle corners at 16  28 and 17  29
\plot 3 30 3 34  /
\plot 8 30 4 34  /
\plot 9 30 9 34  /
\plot 15 30 9.5 34  /
%%%%%%%%%% level 1 %%%%%%%%%%%%%%%%%%%%%%%%%%%%%%%%%%%%%%%%%%%%
\put{$k = 1:$} at  -6 35.5
\put{$\emptyset$} at  3 35.5
\putrectangle corners at 8  35 and 9  36
\plot 3 37 3 41  /
\plot 8 37 4 41  /
%%%%%%%%%% level 0 %%%%%%%%%%%%%%%%%%%%%%%%%%%%%%%%%%%%%%%%%%%%
\put{$k = 0:$} at  -6 42.5
\put{$\emptyset$} at  3 42.5
\endpicture} 
$$
\end{figure}

\begin{rem} {\rm
\begin{itemize}
\item[(a)]  The planar partition diagrams $d \in P_k$ are exactly paired with the paths $(P,Q)$ whose vertices contain only length 0 or 1 partitions.  The planar partition algebra $\CC P_k(n)$  is isomorphic to the Temperley-Lieb algebra $\CC T_{2k}(n)$,  $k \in \frac{1}{2} \ZZ_{>0}$, (see \cite{Jo} or \cite{HR2}) and has the Bratteli diagram shown in Figure \ref{fig:BratteliPlanarPartition}. Our algorithm gives a  combinatorial proof of the following dimension identity for $\CC T_{2k}(n)$,
\begin{equation}
C(2k)  = \sum_{\lambda \in \Lambda_k \atop \ell(\lambda) \le 1}
\left({2k \choose \lfloor k \rfloor - |\lambda|} - {2 k \choose \lfloor k \rfloor - |\lambda|-1}\right)^2, \qquad k \in \frac{1}{2} \ZZ_{>0},
\end{equation}
where $\lfloor ~ \rfloor$ is the floor function and ${2k \choose \lfloor k \rfloor - |\lambda|} - {2 k \choose \lfloor k \rfloor - |\lambda|-1}$  is the dimension of the irreducible $\CC T_{2k}(n)$-module labeled by $\lambda$ (this is also the number of paths to $\lambda$ in the Bratteli diagram of $\CC T_{2k}(n)$).

\item[(b)]  The planar Brauer diagrams $T_k$ span another copy of the  Temperley-Lieb algebra $\CC T_k(n)$ inside $\CC A_k(n)$, and our insertion algorithm when applied to these diagrams gives a different bijective proof of the identity in  (b).

\item[(c)]  In the insertion for planar rook monoid diagrams in $P\!R_k$, we insert or do nothing for steps $j$ with $0 < j \le k$, we delete or do nothing for steps $j + \half$ with $k \le j < 2k$, and we do nothing on all other steps. Thus, the Bratteli diagram, which is shown in Figure \ref{fig:BratteliPlanarPartition}, has levels  $j$ and $j+ \half$ merged and is isomorphic to Pascal's triangle. The irreducible representation associated to the partition of shape $\lambda = (\ell)$ has dimension ${k \choose \ell}$ (see \cite{HH} for the representation theory of $P\!R_k$).  Our algorithm gives an RSK insertion proof of following the well-known binomial identity,
\begin{equation}
{2k \choose k} = \sum_{\ell = 0}^k  {k \choose \ell}^2.
\end{equation}
\end{itemize}}
\end{rem}

Finally, we show that the algorithm has the symmetry property of the usual RSK algorithm for the symmetric group. Namely, if $\pi \to (P,Q)$ for $\pi \in S_k$, then $\pi^{-1} \to (Q,P)$. Using diagrams, the inverse of $\pi$ is simply the diagram of $\pi$ reflected over its horizontal axis.  We define the flip of any partition diagram to be this reflection, so 
for example,
$$
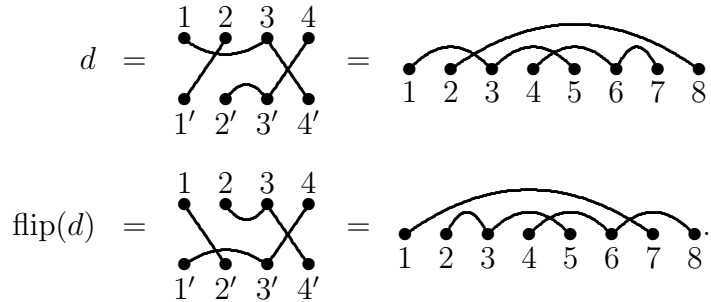
\begin{array}{rcccc}
 d & = &
 {\beginpicture
\setcoordinatesystem units <0.55cm,0.2cm> % sets scale
\setplotarea x from 1 to 4, y from -3 to 3    % sets plot size up
\linethickness=0.5pt                        % sets line thickness
\put{1} at 1 3.5
\put{2} at 2 3.5
\put{3} at 3 3.5
\put{4} at 4 3.5
\put{$1'$} at  1 -3.5
\put{$2'$} at 2 -3.5
\put{$3'$} at 3 -3.5
\put{$4'$} at 4 -3.5
\put{$\bullet$} at 1 -2 \put{$\bullet$} at 1 2
\put{$\bullet$} at 2 -2 \put{$\bullet$} at 2 2
\put{$\bullet$} at 3 -2 \put{$\bullet$} at 3 2
\put{$\bullet$} at 4 -2 \put{$\bullet$} at 4 2
\plot 2 2 1 -2  /
\plot 3 2 4 -2 /
\plot 4 2 3 -2 /
\setquadratic
\plot 1 2 2 1 3 2 /
\plot 2 -2 2.5 -1 3 -2 /
\endpicture} & = &
 {\beginpicture
\setcoordinatesystem units <0.55cm,0.2cm> % sets scale
\setplotarea x from 1 to 4, y from -4 to 3    % sets plot size up
\linethickness=0.5pt                        % sets line thickness
\put{1} at 1 -1.5
\put{2} at 2  -1.5
\put{3} at 3  -1.5
\put{4} at 4  -1.5
\put{$5$} at  5  -1.5
\put{$6$} at 6 -1.5
\put{$7$} at 7  -1.5
\put{$8$} at 8  -1.5
\put{$\bullet$} at 1 0 \put{$\bullet$} at 5 0
\put{$\bullet$} at 2 0 \put{$\bullet$} at 6 0
\put{$\bullet$} at 3 0 \put{$\bullet$} at 7 0
\put{$\bullet$} at 4 0 \put{$\bullet$} at 8 0
\setquadratic
\plot 1 0 2 1.5 3 0 /
\plot 3 0 4 1.5 5 0 /
\plot 4 0 5 1.5 6 0 /
\plot 6 0 6.5 1.5 7 0 /
\plot 2 0 5 3 8 0 /
\endpicture}  \\ \\
\flip(d) & = &
 {\beginpicture
\setcoordinatesystem units <0.55cm,0.2cm> % sets scale
\setplotarea x from 1 to 4, y from -3 to 4    % sets plot size up
\linethickness=0.5pt                        % sets line thickness
\put{1} at 1 3.5
\put{2} at 2 3.5
\put{3} at 3 3.5
\put{4} at 4 3.5
\put{$1'$} at  1 -3.5
\put{$2'$} at 2 -3.5
\put{$3'$} at 3 -3.5
\put{$4'$} at 4 -3.5
\put{$\bullet$} at 1 -2 \put{$\bullet$} at 1 2
\put{$\bullet$} at 2 -2 \put{$\bullet$} at 2 2
\put{$\bullet$} at 3 -2 \put{$\bullet$} at 3 2
\put{$\bullet$} at 4 -2 \put{$\bullet$} at 4 2
\plot 2 -2 1 2  /
\plot 3 -2 4 2 /
\plot 4 -2 3 2 /
\setquadratic
\plot 1 -2 2 -1 3 -2 /
\plot 2 2 2.5 1 3 2 /
\endpicture} 
& = &
 {\beginpicture
\setcoordinatesystem units <0.55cm,0.2cm> % sets scale
\setplotarea x from 1 to 4, y from -3 to 3    % sets plot size up
\linethickness=0.5pt                        % sets line thickness
\put{1} at 1 -1.5
\put{2} at 2  -1.5
\put{3} at 3  -1.5
\put{4} at 4  -1.5
\put{$5$} at  5  -1.5
\put{$6$} at 6 -1.5
\put{$7$} at 7  -1.5
\put{$8$} at 8  -1.5
\put{$\bullet$} at 1 0 \put{$\bullet$} at 5 0
\put{$\bullet$} at 2 0 \put{$\bullet$} at 6 0
\put{$\bullet$} at 3 0 \put{$\bullet$} at 7 0
\put{$\bullet$} at 4 0 \put{$\bullet$} at 8 0
\setquadratic
\plot 8 0 7     1.5 6 0 /
\plot 6 0 5     1.5 4 0 /
\plot 5 0 4     1.5 3 0 /
\plot 3 0 2.5  1.5 2 0 /
\plot 7 0 4     3.0 1 0 /
\endpicture}.
\end{array}
$$
Notice that in the standard representation of $d$, a flip corresponds to a reflection over the vertical line between vertices $k$ and $k+1$.
Our goal is to show that if $d \to (P,Q)$, then $\flip(d) \to (Q,P)$.

Given a partition $d \in A_k$, construct a triangular grid in the integer lattice
$\ZZ \times  \ZZ$ that contains the points in the triangle whose vertices are $(0,0), (2k,0)$ and $(0,2k)$.
Number the columns $1, \ldots, 2k$ from left to right and the rows $1, \ldots, 2k$ from bottom to top.  Place an {\bf X} in the box in column $i$ and row $j$ if and only if, in the one-row diagram of $d$,  vertex $i$ is the left endpoint of edge $j$.  We then label the vertices of the diagram on the bottom row and left column with the empty partition $\emptyset$.  For example, in our diagram $d$ from Example \ref{example:insertionsequence} we have
$$
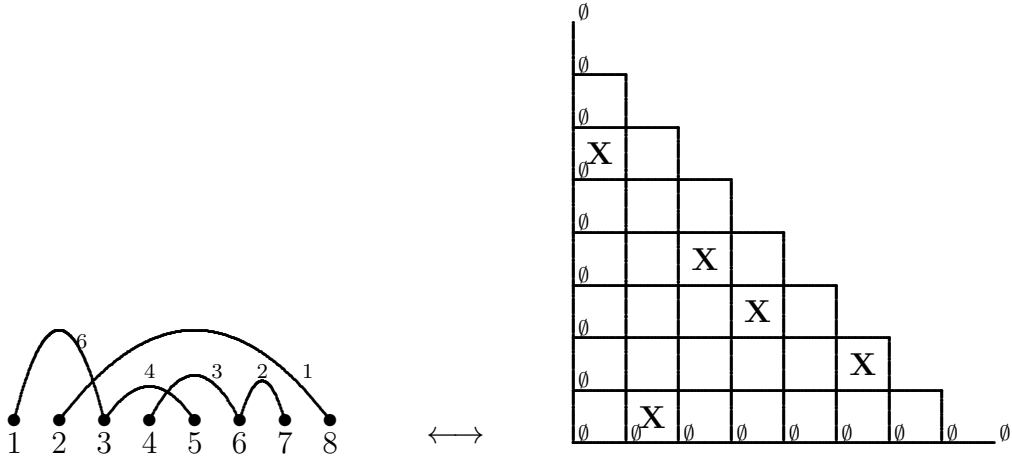
\begin{array}{ccc}
 {\beginpicture
\setcoordinatesystem units <.6cm,0.3cm> % sets scale
\setplotarea x from 0 to 8, y from -2 to 3    % sets plot size up
\linethickness=0.5pt
\put{1} at 1 0
\put{2} at 2 0
\put{3} at 3 0
\put{4} at 4 0
\put{$5$} at 5 0
\put{$6$} at 6 0
\put{$7$} at 7 0
\put{$8$} at 8 0
\put{$\bullet$} at 1 1 \put{$\bullet$} at 7 1
\put{$\bullet$} at 2 1 \put{$\bullet$} at 8 1
\put{$\bullet$} at 3 1 
\put{$\bullet$} at 4 1 
\put{$\bullet$} at 5 1 
\put{$\bullet$} at 6 1
\setquadratic
\plot 1 1   2.0 5.0    3 1 /
\plot 4 1   5 3          6 1 /
\plot 6 1   6.5 2.75   7 1 /
\plot 2 1   5 5          8 1 /
\plot 3 1   4 2.5         5 1 /
\put{$\scriptstyle{1}$} at 7.5 3.2
\put{$\scriptstyle{2}$} at 6.5 3.2
\put{$\scriptstyle{3}$} at 5.5 3.2
\put{$\scriptstyle{4}$} at 4 3.2
\put{$\scriptstyle{6}$} at 2.5 4.5

\endpicture} 
& \qquad \longleftrightarrow\qquad &
{\beginpicture
\setcoordinatesystem units <.7cm,.7cm>    
\setplotarea x from 0 to 8, y from 0 to 8     
\plot 0 0  8 0 /
\plot 0 1  7 1 /
\plot 0 2  6 2 /
\plot 0 3  5 3 /
\plot 0 4  4 4 /
\plot 0 5  3 5 /
\plot 0 6  2 6 /
\plot 0 7    1 7 /

\plot  0 0 0 8 /
\plot  1 0 1 7 /
\plot  2 0 2 6 /
\plot  3 0 3 5 /
\plot  4 0 4 4 /
\plot  5 0 5 3 /
\plot  6 0 6 2 /
\plot  7 0   7 1 /

\put{{\bf X}} at 0.5 5.5
\put{{\bf X}} at 1.5 0.5
\put{{\bf X}} at 2.5 3.5
\put{{\bf X}} at 3.5 2.5
\put{{\bf X}} at 5.5 1.5

\put{$\es$} at .2 .2
\put{$\es$} at 1.2 .2
\put{$\es$} at 2.2 .2
\put{$\es$} at 3.2 .2
\put{$\es$} at 4.2 .2
\put{$\es$} at 5.2 .2
\put{$\es$} at 6.2 .2
\put{$\es$} at 7.2 .2
\put{$\es$} at 8.2 .2

\put{$\es$} at .2 1.2
\put{$\es$} at .2 2.2
\put{$\es$} at .2 3.2
\put{$\es$} at .2 4.2
\put{$\es$} at .2 5.2
\put{$\es$} at .2 6.2
\put{$\es$} at .2 7.2
\put{$\es$} at .2 8.2
 \endpicture}  \\
 \end{array}
 $$
Note that the triangular array completely determines the partition diagram and vice versa.

Now we inductively label the remaining vertices using the local rules of S.\ Fomin (see \cite{Ry1,Ry2} or \cite[7.13]{Sta} and the references there).  If a box is labeled with
$\mu, \nu,$ and $\lambda$ as shown below. Then we add the label $\rho$ according the rules that follow:
$$
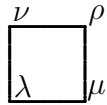
{\beginpicture
\setcoordinatesystem units <.5cm,.5cm>    
\setplotarea x from 0 to 2, y from 0 to 2     
\plot 0 0  2 0 2 2 0 2 0 0 /
\put{$\lambda$} at 0.35 0.4
\put{$\mu$} at 2.3 0.3
\put{$\nu$} at 0.3 2.3
\put{$\rho$} at 2.3 2.3
\endpicture}
$$
\begin{enumerate}
\item[(L1)] If $\mu \not = \nu$, let $\rho = \mu \cup \nu$, i.e., $\rho_i = \max(\mu_i,\nu_i)$ .
\item[(L2)]  If $\mu = \nu$,  $\lambda \subset \mu$, and $\lambda \not= \mu$, then this will automatically imply that $\mu$ can be obtained from $\lambda$ by adding a box to $\lambda_i$. Let $\rho$ be obtained from $\mu$ by adding a box to $\mu_{i+1}$.
 \item[(L3)] If $\mu = \nu = \lambda$, then
if the square does not contain an {\bf X}, let $\rho = \lambda$, and
if the square does contain an {\bf X}, let $\rho$ be obtained from $\lambda$ by adding 1 to $\lambda_1$.
 \end{enumerate}
Using these rules, we can uniquely label every corner, one sep at a time. The resulting diagram is called the {\it growth diagram} $G_d$ for $d$.  The diagram corresponding to the above example is
$$
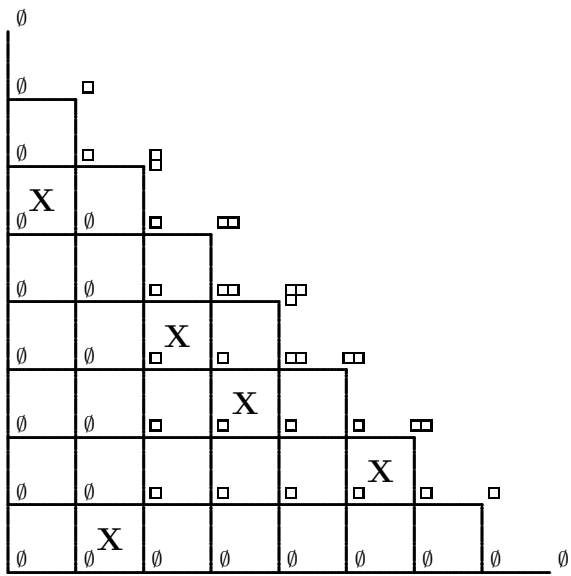
{\beginpicture
\setcoordinatesystem units <.9cm,.9cm>    
\setplotarea x from 0 to 8, y from 0 to 8     
\plot 0 0  8 0 /
\plot 0 1  7 1 /
\plot 0 2  6 2 /
\plot 0 3  5 3 /
\plot 0 4  4 4 /
\plot 0 5  3 5 /
\plot 0 6  2 6 /
\plot 0 7    1 7 /

\plot  0 0 0 8 /
\plot  1 0 1 7 /
\plot  2 0 2 6 /
\plot  3 0 3 5 /
\plot  4 0 4 4 /
\plot  5 0 5 3 /
\plot  6 0 6 2 /
\plot  7 0   7 1 /

\put{{\bf X}} at 0.5 5.5
\put{{\bf X}} at 1.5 0.5
\put{{\bf X}} at 2.5 3.5
\put{{\bf X}} at 3.5 2.5
\put{{\bf X}} at 5.5 1.5

\put{$\es$} at .2 .2
\put{$\es$} at 1.2 .2
\put{$\es$} at 2.2 .2
\put{$\es$} at 3.2 .2
\put{$\es$} at 4.2 .2
\put{$\es$} at 5.2 .2
\put{$\es$} at 6.2 .2
\put{$\es$} at 7.2 .2
\put{$\es$} at 8.2 .2

\put{$\es$} at .2 1.2
\put{$\es$} at .2 2.2
\put{$\es$} at .2 3.2
\put{$\es$} at .2 4.2
\put{$\es$} at .2 5.2
\put{$\es$} at .2 6.2
\put{$\es$} at .2 7.2
\put{$\es$} at .2 8.2

\put{$\es$} at 1.2 1.2
\putrectangle corners at 2.1 1.1 and 2.25 1.25
\putrectangle corners at 3.1 1.1 and 3.25 1.25
\putrectangle corners at 4.1 1.1 and 4.25 1.25
\putrectangle corners at 5.1 1.1 and 5.25 1.25
\putrectangle corners at 6.1 1.1 and 6.25 1.25
\putrectangle corners at 7.1 1.1 and 7.25 1.25

\put{$\es$} at 1.2 2.2
\put{$\es$} at 1.2 3.2
\put{$\es$} at 1.2 4.2
\put{$\es$} at 1.2 5.2
\putrectangle corners at 1.1 6.1 and 1.25 6.25
\putrectangle corners at 1.1 7.1 and 1.25 7.25

\putrectangle corners at 2.1 2.1 and 2.25 2.25
\putrectangle corners at 2.1 3.1 and 2.25 3.25
\putrectangle corners at 2.1 4.1 and 2.25 4.25
\putrectangle corners at 2.1 5.1 and 2.25 5.25
\putrectangle corners at 2.1 6.1 and 2.25 6.25
\putrectangle corners at 2.1 5.95 and 2.25 6.10

\putrectangle corners at 3.1 2.1 and 3.25 2.25
\putrectangle corners at 3.1 3.1 and 3.25 3.25
\putrectangle corners at 4.1 2.1 and 4.25 2.25
\putrectangle corners at 5.1 2.1 and 5.25 2.25

\putrectangle corners at 6.1 2.1 and 6.25 2.25
\putrectangle corners at 5.95 2.1 and 6.1 2.25

\putrectangle corners at 4.1 3.1 and 4.25 3.25
\putrectangle corners at 4.25 3.1 and 4.4 3.25

\putrectangle corners at 5.1 3.1 and 5.25 3.25
\putrectangle corners at 4.95 3.1 and 5.1 3.25

\putrectangle corners at 3.1 4.1 and 3.25 4.25
\putrectangle corners at 3.25 4.1 and 3.4 4.25

\putrectangle corners at 4.1 4.1 and 4.25 4.25
\putrectangle corners at 4.25 4.1 and 4.4 4.25
\putrectangle corners at 4.1 3.95 and 4.25 4.1

\putrectangle corners at 3.1 5.1 and 3.25 5.25
\putrectangle corners at 3.25 5.1 and 3.4 5.25

 \endpicture} .
 $$

Let $P_d$ denote the chain of partitions that follow the staircase path on the diagonal
of $G_d$ from $(0,2k)$ to $(k,k)$ and let $Q_d$ denote the chain of partitions that follow the staircase path on the diagonal of $G_d$ from $(2k,0)$ to $(k,k)$. The pair $(P_d,Q_d)$ represents a pair of vacillating tableaux whose shape is the partition at $(k,k)$.  The staircase path in our growth diagram above is the same as the path we get from insertion in Figure \ref{fig:insertdiag}.

\begin{thm}
\label{thm:growth}
Let $d \in A_k$ with $d \to (P,Q)$. Then $P_d = P$ and $Q_d = Q$.
\end{thm}

\begin{proof}  Turn the diagram $d \in A_k$ into a diagram $d'$ on $4k$ vertices by splitting each vertex $i$ into two vertices labeled by $i-{1 \over 2}$ and $i$. If there is an edge from vertex $j$ to vertex $i$ in $d$, with $j < i$, let $j$ be adjacent to $i-{1 \over 2}$ in $d'$.   If there is an edge from vertex $j$ to vertex $i$ in $d$, with $j > i$, let $i$ be adjacent to $j-{1 \over 2}$ in $d'$.  Thus, in our example we have,
$$
 
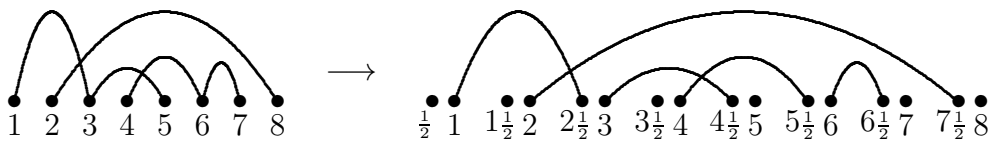
{\beginpicture
\setcoordinatesystem units <.5cm,0.3cm> % sets scale
\setplotarea x from 1 to 8, y from -2 to 3    % sets plot size up
\linethickness=0.5pt
\put{1} at 1 -2
\put{2} at 2 -2
\put{3} at 3 -2
\put{4} at 4 -2
\put{$5$} at 5 -2
\put{$6$} at 6 -2
\put{$7$} at 7 -2
\put{$8$} at 8 -2
\put{$\bullet$} at 1 -1 \put{$\bullet$} at 7 -1
\put{$\bullet$} at 2 -1 \put{$\bullet$} at 8 -1
\put{$\bullet$} at 3 -1 
\put{$\bullet$} at 4 -1 
\put{$\bullet$} at 5 -1 
\put{$\bullet$} at 6 -1
\setquadratic
\plot 1 -1   2.0 3.0 3 -1 /
\plot 4 -1   5 1 6 -1 /
\plot 6 -1   6.5 .75 7 -1 /
\plot 2 -1   5 3 8 -1 /
\plot 3 -1   4 .5 5 -1 /
\endpicture}
\quad\longrightarrow\quad
 {\beginpicture
\setcoordinatesystem units <1cm,0.3cm> % sets scale
\setplotarea x from 1 to 8, y from -2 to 3    % sets plot size up
\linethickness=0.5pt
\put{$\scriptstyle{1 \over 2}$} at .9 -2
\put{1} at 1.3 -2
\put{1$\scriptstyle{1 \over 2}$} at 1.9 -2
\put{2} at 2.3 -2
\put{2$\scriptstyle{1 \over 2}$} at 2.9 -2
\put{3} at 3.3 -2
\put{3$\scriptstyle{1 \over 2}$} at 3.9 -2
\put{4} at 4.3 -2
\put{4$\scriptstyle{1 \over 2}$} at 4.9 -2
\put{$5$} at 5.3 -2
\put{5$\scriptstyle{1 \over 2}$} at 5.9 -2
\put{$6$} at 6.3 -2
\put{6$\scriptstyle{1 \over 2}$} at 6.9 -2
\put{$7$} at 7.3 -2
\put{7$\scriptstyle{1 \over 2}$} at 7.9 -2
\put{$8$} at 8.3 -2
\put{$\bullet$} at 1 -1 \put{$\bullet$} at 7 -1
\put{$\bullet$} at 2 -1 \put{$\bullet$} at 8 -1
\put{$\bullet$} at 3 -1 
\put{$\bullet$} at 4 -1 
\put{$\bullet$} at 5 -1 
\put{$\bullet$} at 6 -1
\put{$\bullet$} at 1.3 -1 \put{$\bullet$} at 7.3 -1
\put{$\bullet$} at 2.3 -1 \put{$\bullet$} at 8.3 -1
\put{$\bullet$} at 3.3 -1 
\put{$\bullet$} at 4.3 -1 
\put{$\bullet$} at 5.3 -1 
\put{$\bullet$} at 6.3 -1

\setquadratic
\plot 1.3 -1   2.15 3.0 3 -1 /
\plot 4.3 -1   5.15 1 6 -1 /
\plot 6.3 -1   6.65 .75 7 -1 /
\plot 2.3 -1   5.15 3 8 -1 /
\plot 3.3 -1   4.15 .5 5 -1 /
\endpicture}
$$
In this way we turn each diagram in $A_k$ into a partial matching in $A_{2k}$.  In his dissertation, T.\ Roby \cite{Ry1} shows that the up-down tableaux coming from growth diagrams are equivalent to those from RSK insertion of Brauer diagrams (matchings).  Dulucq and Sagan \cite{DS} extend this insertion to work for partial matchings (and other generalizations for skew oscillating tableaux) and Roby \cite{Ry2} shows that the growth diagrams also work for the Duluc-Sagan insertion of partial matchings.  Furthermore, the insertion sequences that come from our diagram $d$ are the same as the insertion sequences that come from $d'$ using the methods of  \cite{Sun,DS}.
\end{proof}

A key advantage of the use of growth diagrams is that the symmetry of the algorithm is nearly obvious.  We have that $i$ is the left endpoint of the edge labeled $j$ in $d$  if and only if $j$ is the left endpoint of the edge labeled $i$ in $\flip(d)$.  Thus the growth diagram of $G_d$ is the reflection over the line $y=x$ of the growth diagram of $G_{\flip(d)}$, and so $P_d = Q_{\flip(d)}$ and $Q_d = P_{\flip(d)}$.  Combining this point with the Theorem \ref{thm:growth}  gives

\begin{cor}  If $d \to (P,Q)$ then $\flip(d) \to (Q,P)$.
\label{cor:symmetric}
\end{cor}

We say that $d \in A_k$ is symmetric if $d = \flip(d)$.   The following is an example of a symmetric diagram in $A_6$,
$$
 
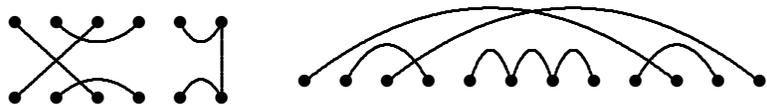
{\beginpicture
\setcoordinatesystem units <0.55cm,0.25cm> % sets scale
\setplotarea x from 0 to 7, y from -2 to 2    % sets plot size up
\linethickness=0.5pt                        % sets line thickness
\put{$\bullet$} at 1 -2 \put{$\bullet$} at 1 2
\put{$\bullet$} at 2 -2 \put{$\bullet$} at 2 2
\put{$\bullet$} at 3 -2 \put{$\bullet$} at 3 2
\put{$\bullet$} at 4 -2 \put{$\bullet$} at 4 2
\put{$\bullet$} at 5 -2 \put{$\bullet$} at 5 2
\put{$\bullet$} at 6 -2 \put{$\bullet$} at 6 2
\plot 1 2 3 -2  /
\plot 3 2 1 -2 /
\plot 6 2 6 -2 /
\setquadratic
\plot 2 2 3 1.0 4 2 /
\plot 2 -2 3 -1.0 4 -2 /
\plot 5 2 5.5 1.0 6 2 /
\plot 5 -2 5.5 -1.0 6 -2 /
\endpicture}
 {\beginpicture
\setcoordinatesystem units <.55cm,0.28cm> % sets scale
\setplotarea x from 0 to 12, y from -1 to 2.5    % sets plot size up
\linethickness=0.5pt
\put{$\bullet$} at 1 -1 \put{$\bullet$} at 7 -1
\put{$\bullet$} at 2 -1 \put{$\bullet$} at 8 -1
\put{$\bullet$} at 3 -1 \put{$\bullet$} at 9 -1
\put{$\bullet$} at 4 -1 \put{$\bullet$} at 10 -1
\put{$\bullet$} at 5 -1 \put{$\bullet$} at 11 -1
\put{$\bullet$} at 6 -1 \put{$\bullet$} at 12 -1
\setquadratic
\plot 1 -1 5.5 2.5 10 -1 /
\plot 3 -1 7.5 2.5 12 -1 /
\plot 2 -1 3 .75 4 -1 /
\plot 9 -1 10 .75 11 -1 /
\plot 5 -1 5.5 .5 6 -1 /
\plot 6 -1 6.5 .5 7 -1 /
\plot 7 -1 7.5 .5 8 -1 /
\endpicture}.
$$
In $S_k$ the symmetric diagrams  are involutions. 

\begin{cor} 
A diagram $d \in A_k$ is symmetric if and only if $d \to (P,P)$.
\end{cor}

\begin{proof}  If $d$ is symmetric, then by Corollary  5.6 we must have $P = Q$.  To prove the converse,
let $P = Q$ and place these vacillating tableaux on the staircase border of a growth diagram. The local rules above are invertible: given $\mu, \nu$ and $\rho$ one can the rules backwards to uniquely find $\lambda$ and determine whether there is an {\bf X} in the box.  Thus, the interior of the growth diagram is uniquely determined. By the symmetry of having $P = Q$ along the staircase, the growth diagram  must have a symmetric interior and a symmetric placement of the {\bf X}s.   This forces $d$ to be symmetric.
\end{proof}

This corollary tells us that the number of symmetric diagrams in $A_k$ is equal to the number of vacillating tableaux of length $2k$, or the number of paths to level $k$ in the Bratteli diagram of $\CC A_k(n)$. Thus,
\begin{equation}\label{eq:model}
\Card(\{ d \in A_k\ | \ d \hbox{ is symmetric } \}) = \sum_{\lambda \in \Lambda_k } m_k^\lambda.
\end{equation}
Furthermore, since our insertion restricts to the diagram subalgebras, (\ref{eq:model}) is true if we replace $A_k$ with any of the  monoids $S_k, B_k, T_k, I_t, R_k, P\!R_k$ and replace $m_k^\lambda$ with the dimension of the appropriate irreducible representation.  In the case of the symmetric group,   (\ref{eq:model}) reduces to the fact that the number of involutions in $S_k$ equals $\sum_{\lambda \vdash k} f^\lambda$.

\end{document}